\input amstex
\documentstyle{amsppt}
\magnification=\magstep1
\pagewidth{6.2truein}
\pageheight{8.3truein}
\NoBlackBoxes
\loadbold
%%%%%%%% macros %%%%%%%%%

\def\A{{\Cal A}}
\def\D{{\Cal D}}
\def\K{{\Cal K}}
\def\L{{\Cal L}}
\def\M{{\Cal M}}
\def\cS{{\Cal S}}
\def\Z{{\Cal Z}}
\def\bC{\text{\bf C}}
\def\bK{\text{\bf K}}
\def\bn{\text{\bf n}}
\def\bR{\text{\bf R}}
\def\bV{\text{\bf V}}
\def\ep{\varepsilon}
\def\bep{\boldsymbol\varepsilon}
\def\dist{\operatorname{dist}}
\def\MIN{\operatorname{MIN}}
\def\sep{\operatorname{SEP}}
\def\csep{\operatorname{CSEP}}
\def\defeq{\mathop{\ \buildrel \text{def}\over =\ }\nolimits}
\def\To{\Rightarrow}
\def\ttT{\tilde{\tilde T}}
%%%%%%%%%%%%%%%%%%%%%%%%%
\topmatter
\title The Complete Separable Extension Property
\endtitle
\author Haskell Rosenthal\endauthor
\abstract 
This work introduces operator space analogues of the Separable Extension 
Property (SEP) for Banach spaces; the Complete Separable Extension Property 
(CSEP) and the Complete Separable Complemention Property (CSCP). 
The results use the technique of a new proof of Sobczyk's Theorem, which also 
yields new results for the SEP in the non-separable situation, e.g., 
$(\oplus_{n=1}^\infty Z_n)_{c_0}$ has the $(2+\ep)$-SEP for all $\ep>0$ 
if $Z_1,Z_2,\ldots$ have the 1-SEP; in particular, 
$c_0 (\ell^\infty)$ has the SEP. 
It is proved that e.g., $c_0(\bR \oplus \bC)$ has the CSEP (where $\bR$, $\bC$ 
denote Row, Column space respectively) as a consequence of the general 
principle: if $Z_1,Z_2,\ldots$ is a uniformly exact sequence of injective 
operator spaces, then $(\oplus_{n=1}^\infty Z_n)_{c_0}$ has the CSEP. 
Similarly, e.g., $\bK_0 \defeq (\oplus_{n=1}^\infty  M_n)_{c_0}$ has the 
CSCP, due to the general principle: $(\oplus_{n=1}^\infty Z_n)_{c_0}$ has 
the CSCP if $Z_1,Z_2,\ldots$ are injective separable operator spaces. 
Further structural results are obtained for these properties, and several 
open problems and conjectures are discussed. 
\vskip8pt
\noindent
{\smc Keywords}: Operator space, local reflexivity, space of compact operators, 
completely complemented subspace, Sobczyk's theorem.
\endabstract 
\endtopmatter

\document
\def\quote#1{\medskip{\parindent=20pt
	\baselineskip=12pt\hsize=5.0truein
	\narrower{\medskip\narrower
	\vbox{\noindent#1}}\smallskip}\smallskip}
\quote{
\line{\hphantom{\S4.} Introduction\dotfill\ \hphantom{2}1\ }
\line{\S1. The Separable Extension Property\dotfill \ \hphantom{2}9}
\line{\S2. The Complete Separable Extension Property \dotfill \ 22}
\line{\S3. The Complete Separable Complementation Property \dotfill\  41}
\line{\S4. Examples of spaces with the CSEP and the CSCP\dotfill\  51}
\line{\hphantom{\S4.} References \dotfill\  55}
}

\baselineskip=18pt				%% Draft mode 

\head Introduction\endhead 

We study here ``quantized'' or ``operator space'' versions of the 
following well known extension property for Banach spaces. 

\proclaim{Definition} 
A Banach space $Z$ is said to have the Separable Extension Property (SEP) 
provided for all separable Banach spaces $Y$, closed linear subspaces $X$, 
and bounded linear operators $T:X\to Z$, there exists a bounded linear 
operator $\tilde T:Y\to Z$ extending $T$. 
\endproclaim

That is, we have the diagram 
$$\matrix 
Y\qquad\quad \cr
\cup\quad\raise.5ex\hbox{$\searrow^{\tilde T}$}\enspace\cr
X\ \mathop{\longrightarrow}\limits_T Z
\endmatrix
\tag 0.1$$

{\it 
If $\lambda\ge1$ is such that $\tilde T$ can always be chosen with 
$\|\tilde T\| \le \lambda\|T\|$, we way $Z$ has the $\lambda$-SEP.} 

In 1941, A.~Sobczyk proved that $c_0$ has the SEP; in fact he showed $c_0$ 
has the 2-SEP, and ``2'' is best possible \cite{S}. 
In 1978, M.~Zippin established the deep converse to this result: 
{\it If $Z$ is a separable infinite-dimensional Banach space with the 
{\rm SEP}, then $Z$ is isomorphic to\/} $c_0$ \cite{Z}. 
These results in a sense ``end'' the study of separable Banach spaces with 
the SEP. 
To the contrary, we show below that ``quantized'' versions of the SEP yield 
a rich ``open-ended'' theory. 
These quantized versions are founded on a new proof for Sobczyk's theorem, 
given in Section~1, which actually yields new information for non-separable 
spaces $Z$ with the SEP. 
For example, we obtain in Corollary~1.4 that $c_0(\ell^\infty)$ 
has the SEP; in fact the $(2+\ep)$-SEP for all $\ep>0$. 
This follows immediately from the new result, Corollary~1.3: 
{\it if $Z_1,Z_2 ,\ldots$ have the {\rm $1$-SEP}, then $(Z_1\oplus Z_2\oplus 
\cdots)_{c_0}$ has the $(2+\ep)$-{\rm SEP}, for all $\ep>0$}. 
A modification  of our argument, due to T.~Oikhberg, actually yields the 
satisfying permanence property: 
{\it if $\lambda\ge1$ and $Z_1,Z_2,\ldots$  have the {\rm $\lambda$-SEP}, 
then $(Z_1\oplus Z_2\oplus \cdots)_{c_0}$ has the 
$(\lambda^2 +\lambda +\ep)$-{\rm SEP} for all $\ep>0$}.
(After circulating the first draft of this paper, I learned that previously 
known results yield that $c_0(\ell^\infty)$ has the 2-SEP --- see 
Remark~3 after Corollary~1.7 below.)

Other proofs of Sobczyk's theorem have been given by 
A.~Pe{\l}czy\'nski \cite{Pe} and W.~Veech \cite{V}. 
We show in Corollary 1.12 below that Veech's argument actually yields the 
isometric result: 
{\it Let  $c_0\subsetneqq Y \subset \ell^\infty$ with $Y$ separable, 
and set $Z= Y/c_0$. 
Then the short exact sequence 
$$0\longrightarrow c_0 \longrightarrow Y\longrightarrow Z\longrightarrow 0$$ 
admits a norm-one lift. 

That is, letting $\pi :Y\to Z$ be the quotient map, there is a norm-one 
operator $L:Z\to Y$ with 
$$\matrix 
Y  \buildrel \pi\over\longrightarrow  Z\cr 
\quad\nwarrow\llap{\hbox{$\scriptstyle L$}}\enspace
\big\uparrow\rlap{\hbox{$\scriptstyle I$}}\cr
\qquad \quad Z\endmatrix$$
Equivalently, $c_0$ is contractively cocomplemented in $Y$; that is, 
there is a linear projection $P$ from $Y$ onto $c_0$ with $\|I-P\| =1$.} 

\noindent 
(This result immediately yields that $c_0$ has the 2-SEP, in virtue of 
the injectivity of $\ell^\infty$.) 
In Theorem~1.9, we use Veech's argument to obtain another generalization 
of Sobczyk's theorem: 
Suppose $X,Y,Z$ are Banach spaces with $Y$ separable and $X\subset Y$, and 
suppose $(T_n)$ is a sequence of operators from $X$ to $Z$ with $T_n\to 0$ 
in the Strong Operator Topology (the SOT). 
If $(T_n)$ admits an extension $(T'_n)$ to $Y$ with $(T'_n)$ relatively compact 
in the SOT, then $(T_n)$ admits an extension $(\tilde T_n)$ to $Y$ with 
$\tilde T_n \to0$ in the SOT. 

{\it That is, setting $T= (T_n)$, $T'= (T'_n)$, $\tilde T= (\tilde T_n)$, 
we have:}

\item{} {\it Hypotheses:} 
$$\matrix Y &\buildrel {T'}\over\longrightarrow  & \ell^\infty (Z)\cr
\cup &&\cup\cr 
X &\buildrel T\over\longrightarrow& c_0(Z)\endmatrix
+ (T'_n) \ \text{{\rm SOT --}  {\it relatively compact.}}$$ 

\item{} {\it Conclusion:} 
$$\matrix 
Y\qquad\qquad\ \  \cr
\cup\quad\raise.5ex\hbox{$\searrow^{\tilde T}$}\qquad\cr
X\ \mathop{\longrightarrow}\limits_T c_0(Z)
\endmatrix\ .$$

The first quantized version of the SEP that we study is the Complete Separable 
Extension Property (CSEP) for operator spaces. 
The definition is obtained by simply inserting ``Completely'' before 
``Separable'' in the definition of the SEP. 
{\it Again, if in the diagram (0.1), $\tilde T$ may always be chosen with 
$\|\tilde T\|_{cb} \le \lambda \|T\|_{cb}$, we say $Z$ has the} 
$\lambda$-CSEP. 

We now briefly recall the following basic concept. 
(For fundamental background and references, see \cite{P} and \cite{Pi}.) 

{\it By an operator space $X$, we mean a Banach space $X$ which is a closed 
linear subspace of $\L(H)$, the bounded linear operator on some Hilbert 
space $H$, endowed with its natural tensor product structure with $\bK$, 
the space of compact operaors on separable infinite dimensional Hilbert space\/}
(which we take as $\ell^2$ for definiteness). 
Thus $\bK \otimes_{op} X$ denotes the closed linear span in $\L(\ell^2 
\otimes_2 H)$ of the operators $K\otimes T$ where $K\in \bK$ and $T\in X$ 
(and $\ell^2\otimes_2H$ is the Hilbert-space tensor product of $\ell^2$ 
and $H$). 
{\it Given operator spaces $X$ and $Y$, a linear operator $T:X\to Y$ is called 
completely bounded if}
$$I\otimes T: \bK \otimes_{op} X\to \bK \otimes_{op}Y$$ 
{\it is bounded; then we set\/} $\|T\|_{cb} = \|I\otimes T\|$. 
It then follows easily that if $X_i$, $Y_i$ are operator spaces and $T_i: 
X_i\to Y_i$ are completely bounded, then $T_1\otimes T_2$ is completely 
bounded, with $\|T_1\otimes T_2\|_{cb} \le \|T_1\|_{cb} \|T_2\|_{cb}$. 
Now many natural Banach space concepts have their operator space versions. 
Thus, operator spaces $X$ and $Y$ are called completely isomorphic if 
there exists an invertible $T:X\to Y$ with $T$ and $T^{-1}$ completely bounded. 
If $\|T\|_{cb} \|T^{-1}\|_{cb} \le\lambda$, we say $X$ and $Y$ are 
$\lambda$-completely isomorphic. 
We then define $d_{cb} (X,Y)$, the completely bounded distance between $X$ 
and $Y$, by 
$d_{cb}(X,Y) = \inf \{ \lambda\ge 1: X$ is $\lambda$-completely 
isomorphic to $Y\}$. 
If $X\subset Y$, with $Y$ an operator space and $X$ a closed linear subspace, 
then $X$ is regarded as an operator subspace of $Y$, via its natural 
structure $\bK\otimes_{op}X\subset \bK \otimes_{op}Y$. 
$X$ is called {\it completely complemented\/} if there is a completely 
bounded projection from $Y$ onto $X$. 
We may then loosely say: 
{\it A separable operator space $X$ has the {\rm CSEP} provided it is 
completely complemented in every separable operator superspace\/}. 
(After the first draft of this paper was completed, it was discovered 
that this ``loose'' statement is actually a theorem, see \cite{OR}.) 

Of course $\bK$ may be identified with a certain Banach space of 
infinite matrices, namely those representing compact operators  on $\ell^2$ 
(with respect to its natural basis). 
For an operator space $X$, $\bK \otimes_{op} X$ may also be visualized 
as a Banach space of infinite matrices, all of whose elements come from $X$. 
We let $M_n$ denote all $n\times n$ matrices of complex scalars, regarded 
as $\L(\ell_n^2)$; we also let $M_{00}$ denote all infinite matrices 
of scalars, with only finitely many non-zero entries. 
Thus we may regard $M_n\subset M_{n+1} \subset \cdots \subset M_{00}
\subset \bK$. 
Now it follows easily that if $P_n :\bK\to M_n$ is the canonical projection, 
then 
$$P_n \otimes I \to I\otimes I\ \text{ in norm, on }\ \bK\otimes_{op} X\ .
\tag 0.2$$ 
For $T:X\to Y$ a bounded linear operator, $n\ge1$, we define $\|T\|_n$ by 
$$\|T\|_n = \|P_n\otimes T\|\ .
\tag 0.3$$ 
(Equivalently, if $I_n=$ Identity in $\L(\ell_n^2)$, 
$\|T\|_n = \|I_n\otimes T\|$.) 
It then follows easily from (0.2) that $T$ is completely bounded iff 
$(\|T\|_n)$ is bounded, and then 
$$\|T\|_{cb} = \sup_n \|T\|_n\ .
\tag 0.4$$ 
(This easy fact is sometimes taken as the definition of complete-boundedness.) 
Visualizing $\bK\otimes_{op} X$ as infinite matrices, we easily then have 
(by the closed graph theorem) that a bounded linear operator $T:X\to Y$ is 
completely bounded exactly when $(Tx_{ij})$ belongs to $\bK \otimes_{op} Y$ for each $(x_{ij})$ in $\bK\otimes_{op} X$; 
of course then $(I\otimes T)$ $(x_{ij}) = (Tx_{ij})$. 

Evidently the concept of an operator space is completely captured by the 
Banach space $\bK \otimes_{op} X$. 
Remarkable axioms of Z.J.~Ruan (cf.\ \cite{ER1}, \cite{Pi}) abstractly 
characterize this tensor product, without reference to the ambient Hilbert 
space. 
Finally, we note that any Banach space $X$ can be regarded as an operator 
space via the so-called MIN structure (where $\|(x_{ij})\|_{\MIN} = 
\sup \{ \|(x^* (x_{ij}))\| : x^* \in X^*$, $\|x^*\|=1\}$ (with 
$\|x^* (x_{ij})\|$ the norm in $\L(H)$). 
Thus formally, Banach space theory is ``subsumed'' by operator space theory. 
However this observation is useless for a Banach space $X$ unless it is 
closely related to $\L(H)$ and its natural subspaces, preduals of such, etc. 
In fact, we can alternatively say that operator space theory is simply a 
special (but very deep!) case of the general theory of tensor products of 
Banach spaces. 

What are some examples of operator spaces with the CSEP? 
Of course $c_0$ has this property; we may ``visualize'' $c_0$ as an operator  
space, by simply identifying it with the space of {\it diagonal\/} matrices 
in $\bK$. 
Similarly, we define $\bR$, the operator Row Space, to be the space of all 
matrices in $\bK$ with entries only in the first row; of course we then 
define $\bC$, the Column Space, as all matrices with entries only in the 
first column. 

It is easily seen that $\bR$ and $\bC$ hve the 1-CSEP. 
We prove (see Corollary after 2.7) that $c_0(\bR\oplus \bC)$ has the 
$(2+\ep)$-CSEP for all $\ep>0$. 
(Throughout, direct sums of operator spaces are taken in the 
$\ell^\infty$-sense.) 
A deep open problem: 
{\it Let $X$ be separable with the\/} CSEP. 
{\it Is $X$ completely isomorphic to a subspace of $c_0 (\bR\oplus \bC)$\/}? 
Of course an affirmative answer would be the direct analogue of Zippin's 
theorem for the CSEP. 
The CSEP structure problem even for subspaces of $c_0(\bR\oplus \bC)$ 
is somewhat involved, however. 
In Section~4, we distinguish 21 (apparently) different infinite-dimensional 
operator subspaces of $c_0(\bR\oplus \bC)$ with the CSEP, represently 
six isomorphically distinct Banach spaces; 
it is conceivable this is the full list (up to complete isomorphism) of all 
infinite-dimensional separable spaces with the CSEP. 

The proof that $c_0(\bR\oplus \bC)$ has the CSEP uses the concept of 
uniformly exact families of operator spaces (see Definition~2.3 below). 
Using Oikhberg's modification of our argument in Section~1 mentioned above, 
we obtain in Corollary~2.5 that {\it if $Z_1,Z_2,\ldots$ are operator spaces 
so that $\{Z_1,Z_2,\ldots\}$ is uniformly exact and the $Z_j$'s all have the 
$\lambda$-{\rm CSEP} for some $\lambda\ge1$, then 
$(Z_1\oplus Z_2\oplus \cdots)_{c_0}$ also has the\/} CSEP. 

We show in Proposition~2.6 that for all $n$, $M_{\infty,n}$ and 
$M_{n,\infty}$ are 1-uniformly exact (where, e.g., $M_{\infty,n}$ denotes 
the $\infty\times n$-matrices in $\bK$). 
Since $M_{\infty,n}$, $M_{n,\infty}$ both have the 1-CSEP, we in fact  obtain 
the following ``almost isometric'' version of the isometric lifting 
property for $c_0$ mentioned above (via Corollary~2.4): 
{\it Fix $n$, and suppose 
$$X = c_0(M_{\infty,n}\oplus M_{n,\infty}) \subsetneqq Y \subset \ell^\infty 
(M_{\infty,n} \oplus M_{n,\infty})$$ 
with $Y$ separable; set $Z= Y/X$. 

\noindent
Then the short exact sequence 
$$0\longrightarrow X\longrightarrow Y \longrightarrow Z\longrightarrow 0$$ 
admits an almost completely contractive lift. 
That is, letting $\pi:Y\to Z$ be the quotient map, then given $\ep>0$, 
there exists an $L:Z\to Y$ with $\|L\|_{cb} <1+\ep$ and} 
$$\matrix 
Y  \buildrel \pi\over\longrightarrow  Z\cr 
\quad\nwarrow\llap{\hbox{$\scriptstyle L$}}\enspace
\big\uparrow\rlap{\hbox{$\scriptstyle I$}}\cr
\qquad \quad Z\endmatrix\qquad\qquad .$$
(This immediately yields that $c_0(M_{\infty,n}\oplus M_{n,\infty})$ 
has the $(2+\ep)$-CSEP for all $\ep>0$, in view of the fact that $\L(H)$ 
is operator-isometrically injective (see \cite{Pi}).) 

We also show that our argument for Theorem~1.9 (via Veech's technique 
\cite{V}) immediately extends to the operator space version (Theorem~2.1a). 
In particular, we obtain for fixed $n$, that 
{\it $c_0(M_n)$ is completely contractively cocomplemented in $Y$ for any 
separable $Y$ with $c_0(M_n) \subsetneqq Y\subset \ell^\infty (M_n)$. 
Thus $c_0(M_n)$ has the\/} 2-CSEP. 

Although the CSEP thus has its isometric aspects, the fundamental constant 
that enters in the CSEP is 2, unless we are already dealing with 
{\it injective\/} operator spaces (i.e., no restrictions on separability 
in the fundamental diagram (0.1) for completely bounded maps). 
Indeed, we show in Proposition~2.22 that if  
{\it $X$ is separable with the {\rm $\lambda$-CSEP} and $\lambda<2$, then 
$X$ is reflexive, whence\/} (via Proposition~2.18) {\it $X$ is 
$\lambda$-injective\/}. 

$\bK$ is often regarded as the quantized version of $c_0$, although 
$\bK_0\defeq (M_1\oplus M_2\oplus \cdots \oplus M_n\oplus \cdots)_{c_0}$ 
is another possible candidate. 
What has become of $\bK$ in our quantization of Sobczyk's theorem? 
E.~Kirchberg establishes in \cite{K} that $\bK$ {\it fails\/} the CSEP. 
In fact, Kirchberg obtains a separable $C^*$-algebra $\A$ and a (two-sided) 
ideal $J\subset\A$ with $J$ *-isomorphic to $\bK$ and $J$ completely 
uncomplemented in $\A$ (and moreover $\A/J$ is an exact $C^*$-algebra). 
Now we note that our positive results hold under a formally weaker 
hypothesis, that of families of operator spaces of finite matrix type (see 
Definition~2.4). 
T.~Oikhberg has established that if conversely $Z_1,Z_2,\ldots$ are 
separable operator spaces and $\{Z_1,Z_2,\ldots\}$ is {\it not\/} of finite 
matrix type, then $(Z_1\oplus Z_2\oplus \cdots)_{c_0}$ fails the CSEP
(see \cite{OR}). 

In particular, $\bK_0$ {\it fails\/} the CSEP (which seems quite surprising 
since $c_0(M_n)$ has the 2-CSEP for all $n$). 
(Actually, it follows from Kirchberg's work \cite{K} and the complete 
isomorphic invariance of exactness for $C^*$-algebras 
(cf.\ \cite{Pi}), that $\bK_0$ fails the CSEP, see Remarks~4.4 \#4, \cite{W}.) 
It turns out that in these counterexamples, the ``culprit'' is the lack of 
local reflexivity of the containing operator space. 

In Section~3, we study a different quantized version of the SEP, the Complete 
Separable Complementation Property (CSCP), which goes as follows: 
{\it A separable operator space $Z$ has the {\rm CSCP} provided every complete 
isomorph of $Z$ is completely complemented in every separable locally 
reflexive operator superspace\/}. 
Equivalently, there exists a completely bounded $\tilde T$ so that the 
diagram (0.1) holds, provided $Y$ is separable locally reflexive and $T$ is 
a complete surjective isomorphism. 
We now indeed obtain that $\bK_0$ {\it has the\/} CSCP. 
It follows from the proof that if e.g., $\A$ is a separable nuclear 
$C^*$-algebra, and $J$ is a $*$-subalgebra $*$-isomorphic to $\bK_0$, then 
$J$ is $(4+\ep)$-completely complemented in $\A$, all $\ep>0$ 
(Corollary~3.2). 

The main result of Section~3, Theorem~3.1, (again via the Oikhberg 
modification mentioned above) goes as follows: 
{\it Let $Z_1,Z_2,\ldots$ be separable operator spaces so that for some 
$\lambda\ge1$, $Z_j$ is $\lambda$-injective for all $j$. 
Then $(Z_1\oplus Z_2\oplus \cdots)_{c_0}$ has the\/} CSCP. 

After the first draft of this paper was completed, it was discovered that 
$\bK$ itself has the CSCP (see \cite{OR}). 
It follows directly from Theorem~3.1 that $(\bigoplus_{n=1}^\infty 
(M_{n,\infty} \oplus M_{\infty,n}))_{c_0}$ has the CSCP. 
Looking at the natural completely complemented subspaces of this space, 
as well as $\bK$ itself, 
it follows that are at least 11 different Banach isomorphism types 
among the separable infinite-dimensional operator spaces with the CSCP, 
and apparently at least 11 different primary such spaces (see Section~4 for 
the relevant definition, and Proposition~4.4 and Conjecture~4.5 for the 
various examples). 
Taking finite direct sums of these, we obtain a finite but apparently 
astronomically large list of separable spaces with the CSCP. 
The deep open question here is thus: 
{\it Does every separable space with the CSCP completely embed in $\bK$?}  
Of course an affirmative answer would yield the ``true'' quantized 
version of Zippin's Theorem \cite{Z}, and also imply it! 
A more accessible problem: 
{\it Classify the infinite-dimensional completely complemented subspaces of 
$\bK$ up to complete isomorphism\/}. 

Some of these results were discovered during the August 1997 Workshop in 
Operator Spaces at Texas A\&M University. 
I wish to thank the many participants with whom I had most useful, 
stimulating discussions. 
Special thanks to Gilles Pisier and Timur Oikhberg for their invaluable 
feedback on the operator-space constructs so crucial in this work.

\head 1. The Separable Extension Property\endhead 

This section is devoted to the pure Banach space category. 

\definition{Definition 1.1} 
A Banach space $Z$ is said to have the {\it Separable Extension Property\/} 
(the SEP) if for all separable Banach spaces $Y$, (closed linear) subspaces 
$X$, and (bounded linear) operators $T:X\to Z$, there exists an operator 
$\tilde T:Y\to Z$ with extending $T$. 
$Z$ is said to have the {\it $\lambda$-SEP\/} provided for all such $X$ 
and $Y$, $\tilde T$ may be chosen with $\|\tilde T\| \le\lambda\|T\|$. 
\enddefinition

We recall also 

\definition{Definition 1.2} 
$Z$ is said to be {\it $\lambda$-injective\/} provided for arbitrary Banach 
spaces $X$ and $Y$ with $X\subset Y$, every operator $T$ from $X$ to $Z$ has 
an extension $\tilde T$ from $Y$ to $Z$ with $\|\tilde T\|\le \lambda \|T\|$. 
\enddefinition 

The results in this section yield new properties of certain non-reflexive 
Banach spaces. 
In the separable space setting, the techniques yield a new proof of 
Sobczyk's theorem that $c_0$ has the SEP \cite{S}, and also yield an intuition 
base for the operator-space discoveries given in the following sections. 
The main ``external'' motivation for the Banach category itself is in 
the non-separable setting, however, because of the profound discovery of 
M.~Zippin \cite{Z}: 
{\it Every infinite dimensional separable Banach space with the $\sep$ is 
isomorphic to $c_0$}. 
We note also that if $Z$ is a finite-dimensional Banach space, then $Z$ has 
the $\lambda$-SEP iff $Z$ is $\lambda$-injective iff $Z$ has the 
$\lambda$-SEP just with respect to finite-dimensional spaces $Y$ and 
subspaces $X$. 
The quantitative characterization of the finite-dimensional 
$\lambda$-injectives  
remains as one of the profound open questions in Banach space theory. 

We also note the classical theorem that 
{\it a Banach space $Z$ is 1-injective iff $Z$ is isomorphic to 
$C(\Omega)$ for some Stonian compact Hausdorff 
space $\Omega$} (real or complex scalars).  
In particular, $\ell^\infty$ is 1-injective and an $n$-dimensional Banach 
space is 1-injective iff it is isometric to $\ell_n^\infty$. 
It also remains a deep open problem whether every infinite-dimensional 
$\lambda$-injective is isomorphic to a 1-injective. 
We note in passing the author's result that $\ell^\infty$ is the 
smallest $\lambda$-injective; i.e., every infinite-dimensional 
$\lambda$-injective contains a subspace isomorphic to $\ell^\infty$ 
\cite{Ro}. 

It can easily be seen (and follows from Proposition 2.1 in the next section) 
that a separable Banach space $Z$ has the SEP iff it has the $\lambda$-SEP 
for some $\lambda\ge1$. 
We do not know if this holds for general Banach spaces $Z$. 
However this is easily seen to be the case if $Z$ is isomorphic to 
$c_0(Z)$ or $\ell^\infty (Z)$. 

Before stating the main result of this section, we give the following 
notation: 
Given $Z_1,Z_2,\ldots$ Banach spaces, $(Z_j)_{c_0}$ denotes the space 
$(Z_1\oplus Z_2\oplus\cdots)_{c_0}$; i.e., the space of all sequences 
$(z_j)$, $z_j\in Z_j$ for all $j$, with $\|z_j\|\to0$; under the natural 
norm $\|(z_j)\| = \sup_j \|z_j\|$. 
Similarly $(Z_j)_{\ell^\infty}$ denotes the space of all bounded sequences 
$(z_j)$ with $z_j\in Z_j$ for all $j$, under the same norm as above. 
Thus $(Z_j)_{c_0}$ is a subspace of $(Z_j)_{\ell^\infty}$. 
Now for a fixed space $Z$, $c_0(Z)$ denotes the space $(Z\oplus Z\oplus
\cdots)_{c_0}$, and similarly $\ell^\infty (Z)=(Z\oplus Z\oplus\cdots)_{c_0}$. 
The reader is thus cautioned that for example, $(\ell_n^2)_{c_0}$ denotes 
the space $(\ell_1^2\oplus \ell_2^2 \oplus \cdots)_{c_0}$, while 
$c_0(\ell_n^2)$ denotes the space $(\ell_n^2 \oplus \ell_n^2
\oplus\cdots)_{c_0}$, for fixed $n$. 

The main result of this section is as follows. 

\proclaim{Theorem 1.1} 
Let $Z_1,Z_2,\ldots$ be 1-injective Banach spaces, $X$, $Y$ be Banach 
spaces with $X\subset Y$ and $Y/X$ separable, and set $Z= (Z_j)_{c_0}$. 
Then for every non-zero operator $T:X\to Z$ and every $\ep>0$, there exists 
a $\tilde T:Y\to Z$ extending $T$ with \newline 
$\|\tilde T\| <(2+\ep)\|T\|$. 
\endproclaim 

Before giving the proof, we give several consequences. 

\proclaim{Corollary 1.2} 
Let $Z$ be as in Theorem 1.2. 
Then $Z$ is $2+\ep$-complemented in every super-space $Y$ with 
$Y/Z$ separable.
\endproclaim 

Recall that if $X\subset Y$, $X$, $Y$ Banach spaces, 
$X$ is called {\it $\lambda$-complemented in\/} $Y$ if there is a (bounded 
linear) projection $P$ mapping $Y$ onto $X$ with $\|P\|\le\lambda$; $Y$ 
is called a super-space of $Z$ if $Z\subset Y$. 
Of course 1.2 follows immediately from 1.1. 

We shall see that the {\it proof\/} of Theorem~1.1 yields 

\proclaim{Corollary 1.3} 
Let $Z_1,Z_2,\ldots$ be Banach spaces with the $1$-$\sep$. 
Then $(Z_j)_{c_0}$ has the $2+\ep$-$\sep$ for all $\ep>0$.
\endproclaim

\remark{Remark} 
T.~Oikhberg has observed that a modification of our argument actually 
yields that if the $Z_j$'s have the $\lambda$-SEP, then $(Z_j)_{c_0}$ has 
the $(\lambda (1+\lambda)+\ep)$-SEP for all $\ep>0$. 
See the Remark following  the proof of Theorem~1.5 below. 
\endremark 

We note in passing that Corollary 1.3 covers certain situations not handled 
by Theorem~1.1. 
For example, let $\Gamma$ be an uncountable set, and let $\ell_c^\infty
(\Gamma)$ denote the space of all scalar-valued functions on $\Gamma$ 
with countable support.  
It is easily seen that $\ell_c^\infty (\Gamma)$ has the 1-SEP; however 
it is known that $\ell_c^\infty (\Gamma)$ is not injective; i.e., not 
$\lambda$-injective for any $\lambda$. 
Corollary~1.3 thus yields that $c_0(\ell_c^\infty (\Gamma))$ {\it has the\/} 
$(2+\ep)$-SEP for all $\ep>0$. 

We note finally the following result, which is perhaps the main 
``external'' motivation for this section.

\proclaim{Corollary 1.4} 
Let $\ep>0$. 
Then $c_0(\ell^\infty)$ has the $(2+\ep)$-$\sep$. 
Moreover $c_0(\ell^\infty)$ is $(2+\ep)$-complemented in every super-space in 
which it is of separable codimension.
\endproclaim

\remark{Remark} 
After circulating the first draft of this paper, I learned that this 
Corollary follows from a known result concerning $M$-ideals in Banach spaces; 
moreover one obtains that ``$\ep$'' may be deleted  in the statment. 
See Remark~3 following the proof of Corollary~1.7 below.
\endremark 

We now pass to the proof of Theorem~1.1. 
We first state a reformulation.
We abuse terminology slightly and say that given Banach spaces $X$, 
$Z_1,Z_2,\ldots$, and operators $T_j:X\to Z_j$ for all $j$, that $(T_n)$ 
tends to zero in the SOT (Strong Operator Topology) provided 
$\|T_nx\|\to 0$ for all $x\in X$. 
(Of course if $Z_1=Z_n$ for all $n$, this is what it means to say $T_n\to0$ 
in the SOT.) 

\proclaim{Theorem 1.5}
Let $Z_1,Z_2,\ldots$, $X$, and $Y$ be Banach spaces satisfying the hypotheses 
of Theorem~1.1. 
Let $(T_j)$ be a sequence of operators with $T_j: Y\to Z_j$ for all $j$, 
so that $T_j\mid X\to 0$ in the {\rm SOT} and $\sup_j\|T_j\|=1$. 
For every $\ep>0$, there exists a sequence $(S_j)$ of operators so that 
for all $j$, 
\roster
\item"1)" $S_j:Y\to Z_j$ 
\item"2)" $X\subset \ker S_j$ 
\item"3)" $\|S_j\| \le 1+ {\ep\over 2}$ 
\endroster
so that $(T_j-S_j)\to 0$ in the {\rm SOT}.
\endproclaim 

We first show that 1.5 $\Rightarrow$ 1.1. 
(The equally easy converse is not needed, and left to the reader. 
Of course the converse motivated the formulation of 1.5.) 

Let $T:X\to Z$ be as in the statement of 1.1; wlg, $\|T\|=1$. 
Let $(T_j)$ be the corresponding sequence so that $Tx= (T_jx)$ 
for all $x\in X$. 
Now for each $j$, since $Z_j$ is 1-injective, choose $T'_j:Y\to Z_j$ 
with $\|T'_j\| = \|T_j\|$ and $T'_j\mid X=T_j$. 
Of course then $\sup_j\|T'_j\| = \sup_j \|T_j\|=1$. 
Let us abuse notation and let $T_j$ denote also the extended operator 
$T'_j$ for all $j$. 
Now the hypotheses of 1.5 hold; choose $(S_j)$ satisfying its conclusion. 
Now define $\tilde T:Y\to (Z_j)_{\ell^\infty}$ by 
$$\tilde T(y) = (T_j-S_j)(y)\ \text{ for all }\ y\in Y\ .
\tag 1.1$$
It follows from the conclusion of 1.5 that $\tilde T$ actually has its 
range in $Z= (Z_j)_{c_0}$, and 2) of 1.5 insures that $\tilde T$ extends $T$; 
of course $\|\tilde T\| <2+\ep$, proving 1.5.\qed 

We next require the following rather surprising 

\proclaim{Lemma 1.6} 
Let $Z_1,Z_2,\ldots$, $X$, $Y$, and  $(T_j)$ be as in 1.5. 
Assume further that $Y/X$ is finite-dimensional and let $F$ be a 
finite-dimensional subspace of $Y$ with $X\oplus F= Y$; let $P$ be the 
projection of $Y$ onto  $F$ with kernel $X$. 
Then 
$$\varlimsup_{n\to\infty} \|T_n P\| \le 1 \ .
\tag 1.2$$
\endproclaim 

\demo{Proof} 
Suppose not. 
By passing to a subsequence, we can suppose wlg there is a $C>1$ so that 
$$\|T_nP\| >C\ \text{ for all }\ n\ .
\tag 1.3$$

So for each $n$, choose $x_n\in X$ and $f_n\in F$ with 
$$\|x_n\oplus f_n\| =1\quad\text{and}\quad \|T_nf_n\| >C\ .
\tag 1.4$$ 
Of course then $\|f_n\|\le \|P\|$ for all $n$; since $F$ is 
finite-dimensional, we can suppose by passing to a further subsequence, that 
$$f_n\to f\ \text{ in norm}
\tag 1.5$$ 
for a certain $f$ in $F$. 
But then $\|T_n (f_n-f)\| \to 0$ as $n\to\infty$, whence by (1.4) and (1.5), 
$$\|T_nf\| >C\ \text{ for all $n$ sufficiently large.}
\tag 1.6$$ 

But also since $x_n\oplus f_n - x_n\oplus f\to 0$ in norm by (1.5), 
$$\lim_{n\to\infty} \|x_n\oplus f\| =1\qquad  \text{ by (1.4).}
\tag 1.7$$ 
Now let $\ep>0$ with $1+\ep <C$, and by (1.7), choose $k$ with 
$$\|x_k\oplus f\| \le 1+\ep\ .
\tag 1.8$$
Then since $\|T_n\|\le 1$ for all $n$, 
$$\|T_n(x_k\oplus f)\| = \|T_n(x_k) + T_n(f)\| \le 1+\ep\ .
\tag 1.9$$
Hence since $\lim_{n\to \infty} \|T_n(x_k)\| =0$ by hypothesis,  
$$\varlimsup_{n\to\infty} \|T_n f\| \le 1+\ep\ .
\tag 1.10$$ 
This of course contradicts (1.6).\qed
\enddemo 

\remark{Remark} 
Lemma 1.6 does not require the hypothesis that the $Z_j$'s be 1-injective, 
and furthermore it immediately yields the conclusion of Theorem~1.5 in case 
$Y/X$ is finite-dimensional. 
Indeed, choose $m$ so that $\|T_nP\| <1+\ep$ for $n\ge m$. 
Let $S_n= 0$ for $n<m$ and $S_n = T_nP$ for $n\ge m$. 
Then $T_n-S_n = T_n (I-P)$ for $n\ge m$, so for $y\in Y$, since 
$(I-P)y\in X$, $\lim_{n\to\infty}\|(T_n-S_n)(y)\| = \|T_n(I-P)(y)\|=0$. 
\endremark 

We thus obtain the following consequence: 

\proclaim{Corollary} 
Let $X\subset Y$ with $Y/X$ finite-dimensional, $(Z_j)$ a given sequence 
of Banach spaces, and $T:Y\to (Z_j)_{\ell^\infty}$ a bounded linear 
operator with $TX\subset (Z_j)_{c_0}$. 
Then given $\ep>0$, there exists $\tilde T:Y\to (Z_j)_{c_0}$ so that 
$\tilde T$ extends $T\mid X$ and $\|\tilde T\|< (2+\ep)\|T\|$.
\endproclaim 

We now pass to the 

\demo{Proof of Theorem 1.5} 
By the above remarks, we may assume that $Y/X$ is infinite-dimensional. 
We may then choose $y_1,y_2,\ldots$ so that $y_1,y_2,\ldots$ are linearly 
independent over $X$ and $Y$ is the closed linear span of $X$ and the $y_n$'s. 
For all $k$, let $F_k = [y_1,\ldots, y_k]$ and $Y_k= X+F_k$. 
Thus $F_k$ is $k$-dimensional and $X\cap F_k=\{0\}$, and 
$$\overline{\bigcup_{k=1}^\infty Y_k} = Y\ .
\tag 1.11$$ 
(Throughout, for any (finite or infinite) sequence $(w_j)$ of elements of a 
Banach space, $[(w_j)]$ denotes the closed linear span of the $w_j$'s.) 

For each $k$, let $P_k:Y_k\to F_k$ be the projection of $Y_k$ onto $F_k$ 
with kernel $X$. 
Let $\ep >0$. 
We shall construct for each $k$, a sequence $(S_n^{(k)})$ of operators 
with the following properties for all $n$; 
$$\gather 
S_n^{(k)} : Y_k \to Z_n \tag 1.12\\
X\subset \ker S_n^{(k)}\tag 1.13\\
\|S_n^{(k)}\| < 1 + {\ep\over 2}\tag 1.14\\
S_n^{(k+1)} \mid Y_k = S_n^{(k)} \tag 1.15
\endgather$$
and also 
$$S_n^{(k)} = T_n P_k \ \text{ for all $n$ sufficiently large.}
\tag 1.16$$

Once $S_n^{(k)}$ has been constructed satisfying (1.12)--(1.16) for all $k$ 
and appropriate $n$'s, for each $n$ let $S_n$ be the unique bounded linear 
operator from $Y$ to $Z_n$ so that 
$$S_n \mid Y_k = S_n^{(k)}\ .
\tag 1.17$$
(As ordered pairs in $Y\times Z_n$, $S_n = \overline{\bigcup_{n=1}^\infty 
S_n^{(k)}}$). 
It follows from (1.14) and (1.15) that $S_n$ is well defined and of course 
$\|S_n\| \le 1+{\ep\over 2} < 1+\ep$; 
(1.13) yields that $X\subset \ker S_n$. 
Of course then $(T_n-S_n)$ is a uniformly bounded sequence of operators; for 
each $k$ and $y\in Y_k$, we have by (1.16) and (1.17) that 
$$(T_n-S_n)(y) = (T_n-T_nP_k) (y) = T_n (I-P_k)(y) 
\tag 1.18$$ 
for all $n$ sufficiently large, whence since $(I-P_k)(y)\in X$, 
$$\|(T_n-S_n)(y)\| \to 0\ .
\tag 1.19$$
Thus $\|(T_n-S_n)(y)\| \to 0$ for all $y\in Y$, since this holds on the dense 
subset $\cup_{k=1}^\infty Y_k$. 
Hence $(S_n)$ satisfies the conclusion of Theorem~1.5. 

We now construct the sequences $(S_n^{(k)})_{n=1}^\infty$ by induction on $k$. 
For convenience, set $S_n^{(0)}=0$ for all $n$. 
Let $k\ge 0$ and suppose $S_n^{(k)}$ has been defined, 
satisfying (1.12)--(1.14) 
for all $n$ and (1.16) for all $n$ sufficiently large. 
By Lemma~1.6, choose $M_{k+1}$ so that for all $n\ge M_{k+1}$ 
$$\|T_nP_{k+1}\| < 1 + {\ep\over 2}
\tag 1.20$$ 
and also (in case $k\ge1$) so that (1.16) holds. 
Now define $S_n^{(k+1)} = T_n P_{k+1}$ for $n\ge  M_{k+1}$. 
For $n<M_{k+1}$, since $Z_n$ is 1-injective, simply choose $S_n^{(k+1)}$ an 
extension of $S_n^{(k)}$ from $Y_k$ to $Y_{k+1}$ with $\|S_n^{(k+1)}\| 
= \|S_n^{(k)}\|$. 

\noindent
(This procedure is also valid in the setting of Corollary~1.3; in this case 
we have that $Y$ will be assumed separable; the assumption that $Z_n$ has 
the 1-SEP again allows us to choose  $S_n^{(k+1)}$ as above.)

We now have that $S_n^{(k+1)}$ satisfies (1.12)--(1.14) and (1.16) (for 
``$k$'' $=k+1$) for all appropriate $n$. 
Finally, we check that (1.15) holds. 
For $n<M_{k+1}$, this is immediate. 
For $n\ge M_{k+1}$, we have since $P_kP_{k+1} = P_k$, that 
$$S_n^{(k)} = T_n P_k = T_n P_{k+1} \mid Y_k = S_n^{(k+1)} \mid Y_k
\tag 1.21$$ 
as desired.\qed
\enddemo 

\remark{Remarks} 
1. A modification of this argument yields that if one instead assumes the 
$Z_j$'s are $\lambda$-injective in Theorem~1.1, then one obtains that 
$\tilde T$ may be chosen as in the conclusion, with $\|\tilde T\| < 
(\lambda(1+\lambda)+\ep)\|T\|$. 
It follows as above that  
{\it if $Z_1,Z_2,\ldots$ have the {\rm $\lambda$-SEP}, then 
$(Z_j)_{c_0}$ has the $(\lambda (1+\lambda)+\ep)$-$\sep$ for all $\ep>0$.}
The modification and these attendant consequences are due to T.~Oikhberg. 
Briefly, assume the $Z_j$'s are $\lambda$-injective and the $T_j$'s as 
in Theorem~1.5, we obtain the $S_j$'s satisfying the conclusion with 
$\|S_j\| < \lambda + \frac{\ep}2$ for all $j$, as follows: 
we construct for each $k$, a sequence $S_n^{(k)}$ of operators so that there 
is an $M_k$ so that for all $n$, 
\roster
\item"(a)" if $n<M_k$ then $S_n^{(k)} : Y\to Z_n$ and $\|S_n^{(k)}\| < 
\lambda +\frac{\ep}2$ 
\item"(b)" if $n\ge M_k$, $S_n^{(k)} :Y_k\to Z_n$, 
$S_n^{(k)} = T_n P_k$, and $\|S_n^{(k)}\| < 1+\frac{\ep}{2\lambda}$ 
\item"(c)" (1.13) and (1.15) hold. 
\endroster
(This modification holds in the complete category also; we give the full 
details in the next section.) 

2. Oikhberg has recently further noted that one may {\it eliminate\/} 
``$\ep>0$'' in the statements of Theorems~1.1 and 1.5, by instead 
constructing $S_n^{(k)}$ in the proof of 1.5 so that $\|S_n^{(k)}\|\le1$ and 
$S_n^{(k)} = (1-2^{-k})T_n P_k$ for $n\ge M_k$ say,  satisfying 1.15 
for $n<M_k$ and 1.12, 1.13 for all $n$. 
The same variation may be used to eliminate ``$\ep>0$'' in the preceding 
remark. 
\endremark 

We note now a further consequence of the proof of Theorem~1.1. 
Given $\lambda \ge 1$ and Banach spaces $X,Y$ with $X\subset Y$, we say 
that $X$ is {\it $\lambda$-cocomplemented in\/} $Y$ if there is a (linear) 
projection $P$ from $Y$ onto $X$ with $\|I-P\| \le\lambda$. 
We say that $X$ is {\it contractively cocomplemented\/} provided $X$ is 
1-cocomplemented; $X$ is {\it almost contractively cocomplemented\/} 
provided $X$ is $(1+\ep)$-cocomplemented for all $\ep>0$. 

\proclaim{Corollary 1.7} 
Let $Z_1,Z_2,\ldots$ be 1-injective Banach spaces, $X= (Z_j)_{c_0}$, and 
let $Y$ 
be a (closed linear) subspace of $(Z_j)_{\ell^\infty}$ with $Y/X$ separable. 
Then $X$ is almost contractively cocomplemented in $Y$.
\endproclaim 

\demo{Proof} 
We easily deduce this from Theorem~1.5. 
Define $T_n:Y\to Z_n$ by $y= (T_1(y),T_2(y),\ldots,T_n(y),\ldots)$ for all 
$y\in Y$. 
($T_n$ is just the restriction of the $n^{th}$ coordinate projection 
on $(Z_j)_{\ell^\infty}$  to $Y$.) 

Let $\ep>0$ and $(S_n)$ be chosen satisfying the conclusion of Theorem~1.5. 
Now defining $P(y) = (T_n(y)- S_n(y))$ for all $y\in Y$, it follows that 
$P$ is a projection from $Y$ onto $X$. 
Indeed, since $T_n-S_n\to 0$ in the SOT, $P$ has its range  in $X$. 
But if $x\in X$, $S_n(x)=0$ for all $n$ and so $x= (T_n(x))_{n=1}^\infty 
= (T_n(x) - S_n(x))_{n=1}^\infty$. 
Of course then $(I-P)(y) = (S_n(y))$ so 
$$\|I-P\| = \sup_n \|S_n\| \le 1 + {\ep\over 2} < 1+\ep\ .
\eqno\qed$$
\enddemo 

\remark{Remarks} 
1. Corollary 1.7 is actually ``stronger'' than Theorem~1.1. 
Indeed, let $X,Y$ and $T:X\to (Z_j)_{c_0}$ be as in the statement of 
Theorem~1.1. 
It follows easily from the 1-injectivity of the $Z_j$'s that also 
$(Z_j)_{\ell^\infty}$ is 1-injective. 
Hence we may choose $T': Y\to (Z_j)_{\ell^\infty}$ extending $T$ so 
that $\|T'\| = \|T\|$. 
Now let $\tilde Y$ denote the closed linear span of $(Z_j)_{c_0}$ and $T'(Y)$. 
Then $\tilde Y /(Z_j)_{c_0}$ is separable, so given $\ep >0$, choose 
$P$ a projection from $\tilde Y$ onto $(Z_j)_{c_0}$ with $\|I-P\| <1+\ep$. 
Of course then $\|P\| < 2+\ep$. 
Now $\tilde T= PT$ is the desired extension of $T$, with 
$\|\tilde T \| < (2+\ep) \|T\|$.\qed

2. Using the modification of the proof of Theorem~1.5 given above, we obtain 
the following generalization: 
{\it Let $Z_1,Z_2,\ldots$ be $\lambda$-injective Banach spaces, and $X$ and 
$Y$ as in Corollary~1.7. 
Then for all $\ep>0$, $X$ is $(\lambda^2+\ep)$-cocomplemented in $Y$\/}.  

3. After the first draft of this paper was completed, it was brought to my 
attention by Bill Johnson that Corollary~1.7 actually follows from a known 
theorem concerning $M$-ideals in Banach spaces, and in fact one obtains 
the stronger conclusion that $X$ {\it is contractively cocomplemented in\/} 
$Y$ ($X,Y$ as in 1.7). 
The theorem, due to T.~Ando \cite{An}, T.~Andersen \cite{A}, 
and later refined by M.~Choi and 
E.~Effros \cite{CE}, yields the following result (see Theorem~2.1, page~59 
of \cite{HWW}). 
Consider a short {\it isometric\/} exact sequence ${0\to X\to Y\to Z\to 0}$. 
{\it Assume that $Z$ is separable and $X$ is an $M$-ideal in $Y$ and an 
$L^1$-predual, then the sequence admits a contractive lift. 
Equivalently, regarding $X\subset Y$, then $X$ is contractively cocomplemented 
in $Y$\/}. 
To obtain Corollary~1.7, we use the known theorem that each $Z_j$ in its 
hypotheses is isometric to $C(\Omega_j)$ for some extremely disconnected 
space $\Omega_j$. 
Making this identification, it then  
follows that $X$ is actually an algebraic closed 
ideal in $(Z_j)_{\ell^\infty}$, which of course may be regarded as a 
$C(\Omega)$-space. 
Thus $X$ is an $M$-ideal and an $L^1$-predual, and the result follows. 
However we note that the generalization of 1.7 given in the previous remark 
does {\it not\/} 
follow from this $M$-ideal result (unless every $\lambda$-injective 
Banach space is isomorphic to a 1-injective, a famous open problem, as  noted 
above). 
\endremark

Lemma 1.6 and the remark following its proof, yield an interesting 
consequence for $(Z_j)_{c_0}$, for general Banach spaces $Z_1,Z_2,\ldots$. 

\definition{Definition 1.3} 
Let $X\subset Y$ be given Banach spaces. 
$X$ is said to be {\it locally complemented\/} 
in $Y$ if there is a $\lambda\ge1$ 
so that 
$$X\text{ is $\lambda$-complemented in $Z$ for all } 
X\subset Z\subset Y \text{ with } Z/X \text{ finite-dimensional.}
\tag 1.22$$
When $(1.22)$ holds, we say $X$ is locally $\lambda$-complemented in $Z$.
\enddefinition

\proclaim{Corollary 1.8} 
Let $Z_1,Z_2,\ldots$ be arbitrary Banach spaces. 
Then $(Z_j)_{c_0}$ is locally almost contractively cocomplemented in 
$(Z_j)_{\ell^\infty}$. 
\endproclaim 

The statement means that $(Z_j)_{c_0}$ is almost contractively cocomplemented 
in $Y$ for all $Y$ with $(Z_j)_{c_0}\subset Y\subset (Z_j)_{\ell^\infty}$ 
and $Y/(Z_j)_{c_0}$ finite-dimensional. 
{\it Thus $(Z_j)_{c_0}$ is locally $(2+\ep)$-complemented in 
$(Z_j)_{\ell^\infty}$, for all $\ep>0$}. 
The proof of 1.8 follows immediately from the Remark following the proof 
of Lemma~1.6, and the argument for Corollary~1.7. 

\remark{Remark} 
W.~Johnson and T.~Oikhberg have obtained a stronger result 
when the $Z_j$'s are separable 
\cite{JO}; see the remark following Corollary~1.12 below. 
\endremark 

As noted above, our argument yields a new proof of Sobczyk's Theorem that 
$c_0$ has the SEP, but we pay an ``$\ep$'' price, for in fact $c_0$ has 
the 2-SEP by \cite{S}, (and ``2'' is best possible here). 
We recapture this result through the following extension theorem, 
whose proof uses a technique due to W.~Veech \cite{V}. 

\proclaim{Theorem 1.9} 
Let $Z$ be an arbitrary Banach space, $X\subset Y$ separable Banach spaces, 
and $T:Y\to \ell^\infty (Z)$ a bounded linear operator so that 
\roster 
\item"(i)" $T(X) \subset c_0(Z)$
\endroster
and
\roster
\item"(ii)" $(T_n)$ is relatively compact in the {\rm SOT}, where 
$T(y) = (T_n(y))$ for all $y$.
\endroster
Then there exists an operator $\tilde T:Y\to c_0(Z)$ extending $T\mid X$ 
with $\|\tilde T\| \le 2\|T\|$. 
\endproclaim 

\remark{Remark}
The hypotheses hold for any finite-dimensional $Z$. 
This easily yields the fact that $c_0$ has the 2-SEP; we give the detailed 
proof in Corollary~2.14 below. 
\endremark

Theorem 1.9 is a consequence of the following two simple lemmas. 

\proclaim{Lemma 1.10} 
There exists a norm $|\cdot|$ on $\L (Y,Z)$ so that letting $\M=(\L(Y,Z),
|\cdot|)$, then the $\M$-topology coincides with the {\rm SOT} topology 
on bounded subsets of $\L(Y,Z)$. 
\endproclaim 

\demo{Proof} 
Let $d_1,d_2,\ldots$ be a countable dense subset of the unit ball of $Y$ and 
define $|\cdot|$ by 
$$|S| = \sum_{n=1}^\infty { \|S(d_n)\|\over 2^n}\ .
\tag 1.23$$
It is easily verified that $\M$ is a normed linear space. 
Moreover, if $(T_n)$ is a bounded sequence in $\L(Y,Z)$ and $T\in \L(Y,Z)$, 
then $T_n\to T$ SOT iff $T_n(d_j)\to T(d_j)$ all $j$ iff 
$|T_n-T| \to 0$.\qed
\enddemo 

\proclaim{Lemma 1.11} 
Let $(\M,\rho)$ be a metric space, $\cS$ a closed subset of $\M$, and $(T_n)$ 
a sequence in $\M$ so that $\{T_1,T_2,\ldots\}$ is relatively compact and 
all cluster points of $(T_n)$ lie in $\cS$. 
There exists a sequence $(S_n)$ of points in $\cS$ so that 
$$\rho (T_n,S_n)\to 0\ \text{ as }\ n\to\infty\ .
\tag 1.24$$
\endproclaim

\demo{Proof} 
Define $d_n$ by 
$$d_n = \dist (T_n,\cS) \defeq \inf \{\rho (T_n,S):S\in \cS\}\ .
\tag 1.25$$
Of course (1.24) is simply the assertion that $d_n\to 0$ as $n\to \infty$. 
Were this false, by passing to a subsequence if necessary, we can assume 
wlg there is a $d>0$ so that 
$$d_n \ge d\ \text{ for all }\ n\ .
\tag 1.26$$
Choose $n_1 < n_2<\cdots$ so that $(T_{n_i})$ converges, to $S$ say. 
By hypothesis, $S\in \cS$. 
Now $\rho (T_{n_i},S)\to 0$ as $i\to \infty$, yet $\rho (T_{n_i},S)\ge 
\dist (T_{n_i},\cS)\ge d_{n_i} \ge d>0$ for all $i$, 
a contradiction.\qed
\enddemo 

\demo{Proof of Theorem 1.9} 
We may assume wlg that $\|T\|=1$. 
Let $\M$ be as in Lemma~1.10, and define $\cS$ by 
$$\cS = \{ S\in \L(Y,Z) :\|S\| \le 1\text{ and } X\subset \ker S\}\ .
\tag 1.27$$
Now the hypotheses imply that all SOT-cluster points of $(T_n)$ lie in $\cS$. 
Indeed, if $T_{n_i}\to T$ in the SOT, then for $x\in X$, 
$T_{n_i}(x) \to T(x)$, but  $T_{n_i}(x)\to 0$ in norm; of course 
$\|T\|\le 1$ since $\|T_j\| \le 1$ for all $j$. 
Thus the hypotheses of Lemma~1.11 apply (where of course $\M$ is endowed 
with the standard metric, $\rho (x,y) = |x-y|$). 
Hence we may choose a sequence $(S_n)$ in $\cS$ so that 
$|T_n-S_n| \to 0$ as $n\to\infty$; i.e., by Lemma~1.10, 
$$T_n-S_n \to 0\ \text{ in the SOT}. 
\tag 1.28$$ 
Now define $\tilde T$ by $\tilde T(y) = (T_n-S_n)(y)$ for all $y\in Y$. 
Then $\tilde T$ is the desired extension of $T\mid X$.\qed
\enddemo 

We may now deduce the following rather surprising consequence of this proof, 
analogous to Corollary~1.7. 

\proclaim{Corollary 1.12}
Let $Z$ be a finite-dimensional Banach space and $Y$ a separable subspace 
of $\ell^\infty (Z)$ containing $c_0(Z)$. 
Then $c_0(Z)$ is contractively cocomplemented in $Y$. 
\endproclaim 

\demo{Proof} 
Of course we assume $Y\ne c_0(Z)$. 
Let $T$ denote the identity injection of $Y$ into $\ell^\infty (Z)$, and let 
$y= (T_1(y),T_2(y),\ldots,T_n(y),\ldots)$ for all $y\in Y$. 
Since $Z$ is finite-dimensional and of course $T_n:Y\to Z$ with 
$\|T_n\| \le 1$ for all $n$, $(T_n)$ is relatively compact in the SOT on 
$\L(Y,Z)$. 
Now choose $(S_n)$ as in the proof of Theorem~1.9. 
It follows that defining $P:Y\to c_0(Z)$ by 
$$P(y) = (T_n)(y) - S_n(y)) \ \text{ for all }\ y\in Y\ ,$$ 
then $P$ is a linear projection from $Y$ onto $c_0(Z)$, and 
$(I-P)(y) = (S_n(y))$ for all $y$, whence $\|I-P\| =\sup \|S_n\| =1$.\qed
\enddemo 

\remark{Remark} 
This result has been obtained in \cite{JO}, using a different argument. 
The authors of \cite{JO} also study the family of separable Banach spaces $Z$ 
so that $Z$ is complemented in $Y$ for all separable $Y$ with 
$c_0(Z) \subset Y \subset \ell^\infty (Z)$, obtaining quite nice  
results, including the fact that there exists a sequence $(E_n)$ 
of finite dimensional Banach spaces so that $Z= (E_n)_{c_0}$ {\it fails\/} 
this property. 
It is also proved in \cite{JO} that for $Z$ separable, $c_0(Z)$ is locally 
contractively cocomplemented in $\ell^\infty (Z)$, thus removing the 
``almost'' from Corollary~1.8 above, when the $Z_j$'s are separable. 
\endremark 

As noted in the introduction, we immediately obtain the following 
isometric property for $c_0$ itself. 

\proclaim{Corollary 1.13} 
$c_0$ is contractively cocomplemented in any separable superspace 
which lies in $\ell^\infty$. 
\endproclaim 

In turn, this yields Sobczyk's Theorem. 

\proclaim{Corollary 1.14} 
$c_0$ has the $2$-$\sep$. 
\endproclaim 

\demo{Proof} 
Let $X\subset Y$ be separable Banach spaces and let $T:X\to c_0$ be a given 
bounded linear operator. 
Define the sequence $(f_n)$ in $X^*$ by $Tx= (f_1(x),f_2(x),\ldots)$ for 
all $x$. 
Now fixing $n$, $\|f_n\| \le \|T\|$; let $f'_n$ be a Hahn-Banach extension 
of $f_n$ to $Y$. 
Let $T' : Y\to \ell^\infty$ be defined by $T'(y) = (f'_1(y),f'_2(y),\ldots)$. 
Then of course $\|T'\| = \|T\|$. 
Let $\tilde Y$ denote the closed linear span of $Y$ and $c_0$, and choose 
$P$ a projection from $\tilde Y$ onto $c_0$ with $\|I-P\|=1$. 
Then $\tilde T= PT'$ yields the desired extension of $T$ to $Y$ with 
$\|\tilde T\| \le 2\|T\|$.\qed
\enddemo 

\remark{Remark} 
Corollary 1.13 also follows directly from the $M$-ideal theorem cited in 
Remark~3 following the proof of Corollary~1.7. 
\endremark 

\head 2. The Complete Separable Extension Property\endhead 

As noted in the Introduction, most of the results of the previous section 
follow from their operator space versions given here. 
However the {\it techniques\/} of proof come from the arguments in Section~1, 
so we have chosen to present the Banach category first, in the interest 
of clarity. 

We first recall the definition given in the Introduction. 

\definition{Definition 2.1} 
An operator space $Z$ is said to have the Complete Separable Extension 
Property (the {\rm CSEP}) if for all separable operator spaces $Y$, 
subspaces $X$, and completely bounded operators $T:X\to Z$, there exists 
a completely bounded operator $\tilde T: Y\to Z$ extending $T$. 
$Z$ is said to have the $\lambda$-CSEP provided for all such $X$ and $Y$, 
$\tilde T$ may be chosen with $\|\tilde T\|_{cb} \le\lambda \|T\|_{cb}$. 
\enddefinition 

Again, we have the operator space analogue of injectivity in Banach spaces. 

\definition{Definition 2.2} 
An operator space $Z$ is said to be {\it isomorphically injective\/}
provided for 
arbitrary  operator spaces $X$ and $Y$ with $X\subset Y$, every completely 
bounded  map $T:X\to Z$ admits a completely bounded extension 
$\tilde T :Y\to Z$. 
$Z$ is called {\it $\lambda$-injective\/} 
if the extension $\tilde T$ may always be 
chosen with $\|\tilde T\|_{cb} \le \lambda \|T\|_{cb}$. 
Finally, $Z$ is called {\it isometrically injective\/} if it is $1$-injective.
\enddefinition 

It is a basic theorem in operator space theory that $\L(H)$ is 
isometrically injective for all Hilbert spaces $H$. 
The theorem was proved for the fundamental case of completely positive 
maps and self-adjoint operator spaces in the domain by W.B.~Arveson 
\cite{Arv}, and later in general by V.~Paulsen (cf.\ \cite{P}) and 
G.~Wittstock \cite{Wi}. 
See also \cite{Pi} for a proof from the abstract operator-space viewpoint. 
It follows easily that an operator space $X$ is isomorphically injective 
provided it is completely complemented in some complete isometric embedding 
$\tilde X$ into $\L(H)$; moreover if $P :\L(H) \to \tilde X$ is a completely 
bounded projection, then $X$ is $\lambda$-injective where $\lambda = 
\|P\|_{cb}$. 
We will mainly be concerned with isometrically injective operator spaces here. 
Unlike the Banach space category, there are separable infinite-dimensional 
examples. 
A complete classification of these has been given by G.~Robertson \cite{Ro} 
(see Section~4 below). 
See also work of Z.J.~Ruan giving a characterization of the 1-injectives 
as ``corners'' of injective $C^*$-algebras \cite{Ru}. 
Finally, we note that the 1-CSEP is studied for $C^*$-algebras by 
R.R.~Smith and D.P.~Williams \cite{SW}. 

The 1-injectivity of $\L(H)$ easily yields the following result. 

\proclaim{Proposition 2.1} 
Let $X$ be a separable operator space. 
The following are equivalent 
\roster 
\item"1)" $X$ has the {\rm CSEP}.
\item"2)" $X$ is completely complemented in every separable operator space 
$Y$ with $X\subset Y$. 
\item"3)" There is a $\lambda \ge 1$ so that $X$ is $\lambda$-completely 
complemented in every separable operator space $Y$ with $Y\supset X$. 
\endroster
Moreover if $X$ satisfies {\rm 3)}, $X$ has the {\rm $\lambda$-CSEP}. 
\endproclaim 

\demo{Proof} 

$1) \To 2)$. Trivial. 

$2) \To 3)$. We may assume $X\subset \L(H)$ (with $H$ separable 
infinite-dimensional Hilbert space). 
We shall prove 
\roster
\item"$3')$"
{\it there is a $\lambda\ge1$ so that $X$ is completely 
$\lambda$-complemented in separable super space of $X$ contained in $\L(H)$.}
\endroster 
Were this false, we could choose $Y_1,Y_2,\ldots$ separable operator 
subspaces of $\L(H)$, so that for all $n$, $X\subset Y_n$ but $X$ is not 
completely $n$-complemented in $Y_n$. 
Then letting $Y= [Y_1,Y_2,\ldots]$, $Y$ is separable, $X\subset Y\subset 
\L(H)$, but $X$ is not-completely complemented in $Y$, a contradiction. 

Now suppose $X$ satisfies $3')$. 
Let then $Y\subset Z$ be separable operator spaces and $T:Y\to X$ a 
completely bounded operator. 
Choose $\tilde T: Z\to \L(H)$ extending $T$, with $\|\tilde T\|_{cb} = 
\|T\|_{cb}$ (by the fundamental theorem cited above). 
Now letting $E= \overline{X+TY}$, $E$ is separable, 
so choose $P:E\to X$ a completely 
bounded projection with $\|P\|_{cb} \le \lambda$. 
Then $T' \defeq P\tilde T$ is an extension of $T$ to $Z$, and $\|T'\|_{cb} 
\le \|P\|_{cb} \|\tilde T\|_{cb} \le \lambda \|T\|_{cb}$, proving that $X$ 
has the $\lambda$-CSEP (so of course 3) holds).\qed
\enddemo 

\remark{Remark} 
As in the Banach space category, we do not know if (non-separable) 
operator spaces $X$ with the CSEP have the $\lambda$-CSEP for some  
$\lambda <\infty$. 
Again, this is easily seen to be true if $X$ is completely isomorphic 
to $c_0(X)$ or $\ell^\infty (X)$. 
\endremark 

We next pass to a rather strong condition on operator spaces, which we 
will use to produce examples of spaces with the CSEP. 
 
\definition{Definition 2.3} 
A family $\Z$ of operator spaces is called {\it uniformly exact\/} if 
there is a $C\ge 1$ and a function $\bn : N\to N$ so that for all 
$Z\in \Z$ and all $k$-dimensional subspaces $F$ of $Z$, there exists 
a $G\subset M_{\bn (k)}$ with 
$$d_{cb} (F,G) \le C\ .
\tag 2.1$$
In case $C$ works, we say $\Z$ is {\it $C$-uniformly exact\/}. 
In case $\bn$ works, we say $\bn$ {\it is a uniformity function for\/} $\Z$. 
We say an operator space $Z$ is {\it uniformly exact\/}  (resp. 
$C$-uniformly exact) in case $\Z = \{Z\}$ has the corresponding property. 
\enddefinition 

It follows that an operator space $Z$ is $C$-exact as defined in \cite{Pi}  
precisely when every finite dimensional subspace of $Z$ is $C+\ep$-uniformly 
exact for every $\ep>0$. 
If $X$ is a Banach space endowed with the MIN operator space structure, 
then $X$ is $(1+\ep)$-uniformly exact for every $\ep>0$. 
We may now state the main result of this section, (which yields 
Theorem~1.1 in view of the last comment above). 

\proclaim{Theorem 2.2} 
Let $Z_1,Z_2,\ldots$ be $\lambda$-injective operator spaces so that 
$\{Z_1,Z_2,\ldots\}$ is $C$-uniformly exact for some $C\ge 1$, and set 
$Z = (Z_j)_{c_0}$. 
Let $X\subset Y$ be operator spaces with $X\subset Y$ and $Y/X$ separable. 
Then for every non-zero completely bounded operator $T:X\to Z$ and 
every $\ep >0$, there exists a completely bounded operator 
$\tilde T :Y\to Z$ extending $T$ with 
$\|\tilde T\|_{cb} < (C\lambda^2+ \lambda +\ep) \|T\|_{cb}$. 
\endproclaim 

\remark{Remark} 
We had originally obtained this result for $\lambda=1$. 
This more general result follows via the proof-modification due to Oikhberg, 
mentioned in Section~1. 
\endremark 

We again give several consequences before passing to the proof. 

\proclaim{Corollary 2.3} 
Let $Z$ be as in Theorem 2.2 and let $\ep>0$. 
Then $Z$ is completely $(C\lambda^2 +\lambda +\ep)$-complemented in 
every operator superspace $Y$ with $Y/Z$ separable.
\endproclaim 

Of course this corollary follows immediately from Theorem~2.1. 
Inserting ``completely'' before the ``cocomplemented'' definition given 
preceding Corollary~1.7, we again discover the following consequence 
of the proof of Theorem~2.2. 

\proclaim{Corollary 2.4} 
Let $Z_1,Z_2,\ldots$ satisfy the hypotheses of Theorem~2.1, $Z=(Z_j)_{c_0}$,  
and let $Y$ be an operator space with $(Z_j)_{c_0}\subset Y\subset 
(Z_j)_{\ell^\infty}$ and $Y/Z$ separable and let $\ep>0$. 
Then $Z$ is $(C\lambda^2 +\ep)$-completely cocomplemented in $Y$. 
\endproclaim 

We may thus conclude, 
{\it if $\{Z_1,Z_2,\ldots\}$ is $(1+\ep)$-uniformly exact for every 
$\ep>0$ and the $Z_j$'s are $1$-injective,  
then $Z$ is almost completely contractively cocomplemented in $Y$}. 

Again, we have the following analogue of Corollary~1.3 (which follows 
immediately from Theorem~2.1 if the $Z_j$'s are 1-injective). 

\proclaim{Corollary 2.5} 
If $Z_1,Z_2,\ldots$ are operator-spaces with the $\lambda$-$\csep$ and 
$\{Z_1,Z_2,\ldots\}$ is $C$-uniformly exact, then $(Z_j)_{c_0}$ has the 
$(C\lambda^2 +\lambda +\ep)$-$\csep$ for all $\ep>0$.
\endproclaim 

The next result, combined with Theorem 2.2, yields our examples of 
separable spaces with the CSEP. 

\proclaim{Proposition 2.6} 
For all $n$, $M_{\infty,n}$ and $M_{n,\infty}$ are 1-uniformly exact, with 
uniformity function $\bn (k) = k\cdot n$. 
\endproclaim 

\remark{Remark} 
Every isometrically injective separable operator space $Z$ is completely 
isometric to a subspace of $M_{\infty,n} \oplus M_{n,\infty}$ for some $n$. 
It thus follows from Corollary~2.4 that for such $Z$'s that $c_0(Z)$ 
{\it is almost completely contractively cocomplemented in $Y$ for all 
separable $Y$ with\/} $c_0 (Z) \subset Y \subset \ell^\infty (Z)$. 
\endremark 

The following corollary gives the main ``separable'' motivation for 
Theorem~2.2. 

\proclaim{Corollary 2.7} 
For all $n$, $c_0$ $(M_{\infty,n} \oplus M_{n,\infty})$ has the 
$(2+\ep)$-$\csep$ for all $\ep>0$. 
\endproclaim 

\remark{Remark} 
As noted in the introduction, T.~Oikhberg has shown that 
$(M_n)_{c_0}$ fails the CSEP. 
In view of this, it seems rather surprising that the above family of spaces 
has the CSEP with a good uniform constant. 
\endremark 

Of course 2.7 yields the immediate 

\proclaim{Corollary} 
$c_0(\bR \oplus \bC)$ has the $(2+\ep)$-$\csep$ for all $\ep>0$. 
\endproclaim 

I do not know if the ``$\ep>0''$ can be removed 
from this statement, but I guess the answer is no. 

The fundamental open problem for the {\it characterization\/} of separable 
spaces with the CSEP goes as follows: 

\proclaim{Problem} 
Let $X$ be a separable operator space with the $\csep$. 
Is $X$ completely isomorphic to a completely contractively complemented 
subspace of $c_0 (\bR\oplus \bC)$?
\endproclaim 

A more ``refined'' version of this problem is given in Section~4 below. 

It seems very likely that the fundamental problem reduces to the 

\proclaim{Embedding Problem} 
Let $X$ be a separable operator space with the $\csep$. 
Is $X$ completely isomorphic to a subspace of $c_0(\bR\oplus \bC)$? 
\endproclaim 

In turn, Corollary 2.4 leads us to the following 

\proclaim{Quantitative Embedding Problem} 
Let $\lambda \ge1$, and let $X$ be  separable with the $\lambda$-$\csep$. 
Is there a $\beta$ depending only on $\lambda$, and an $n$ 
(depending on $X$), so that $d_{cb} (Y,X) \le\beta$ for some subspace $Y$ 
of $c_0(M_{\infty,n} \oplus M_{n,\infty})$? 
\endproclaim 

Before dealing with the proofs of the results stated above, we give some 
information concerning the relationship between the CSEP and isomorphic  
injectivity for operator spaces. 
The following result follows quickly from known theorems. 

\proclaim{Proposition 2.8} 
Let $X$ be a non-reflexive  operator space. 
If $X$ is completely isomorphic to a completely complemented subspace of some 
$C^*$-algebra, then $X$ contains a subspace completely isomorphic to $c_0$. 
If moreover $X$ is completely isomorphic to a completely complemented subspace 
of some von-Neumann algebra, then $X$ contains a subspace completely 
isomorphic to $\ell^\infty$. 
\endproclaim 

\remark{Comment} 
This result also holds if one deletes the term ``completely'' from all 
occurrences in its statement.
\endremark 

\demo{Proof} 
Suppose first without loss of generality that $X\subset \A$, $\A$ a 
$C^*$-algebra, and $P:\A\to X$ is a completely bounded projection onto $\A$. 
Since $P$ is non-weakly compact, a result of H.~Pfister \cite{Pf}  
yields there exists 
a commutative $C^*$-subalgebra $\tilde{\A}$ of $\A$ with $P\mid \tilde{\A}$ 
non-weakly compact. 
By uniqueness of the operator-space structure for $C^*$-algebras, it follows 
that $\tilde{\A}$ has MIN as its inherited operator space structure. 
By a result of A.~Pe{\l}czy\'nski \cite{Pe}, there exists a subspace $E$ of 
$\tilde{\A}$ with $E$ isomorphic (and hence completely isomorphic) to $c_0$, 
so that $P\mid E$ is a Banach isomorphism. 
Since $P$ is completely bounded and $E$ has the MIN structure, $P|E$ 
is in fact a complete isomorphism onto its range, proving the first assertion. 
Now if $\A$ is a von-Neumann algebra, then of course $\A$ is also a dual 
Banach space, and in fact the canonical projection $\Pi : \A^{**}\to \A$ 
is then completely bounded. 
Now of course $E^{**}$ is completely isomorphic to $\ell^\infty$; 
regarding $E^{**} = E^{\bot\bot} \subset \A^{**}$, the operator $T$ 
defined by $T= P\Pi| E^{**}$ is then a non-weakly compact completely bounded 
operator into $X$. 
By a result of the author \cite{Ro}, it follows that there is a subspace $Z$ of 
$E^{**}$ with $Z$ isomorphic to $\ell^\infty$ and $T|Z$ a Banach isomorphism. 
Again, $E^{**}$ has the MIN structure, hence so does $Z$, so as before, $T|Z$ 
is a complete isomorphism of $Z$ 
onto its range.\qed
\enddemo 

The following result is now immediate. 

\proclaim{Corollary 2.9} 
Let $X$ be a non-reflexive operator space. 
If $X$ has the $\csep$, $X$ contains a subspace completely isomorphic to $c_0$. 
If $X$ is isomorphically injective, $X$ contains a subspace completely 
isomorphic to $\ell^\infty$. 
\endproclaim 

The next result shows that {\it reflexive\/} separable operator spaces 
are isomorphically injective provided they have the CSEP. 

\proclaim{Proposition 2.10} 
Let $X$ be a reflexive separable operator space with the $\lambda$-$\csep$. 
Then $X$ is a $\lambda$-injective operator space.
\endproclaim 

\demo{Proof}
Assume $X\subset \L(H)$ (with $H$ separable infinite-dimensional Hilbert space) 
and let $\cS$ be the family of separable subspaces $Y$ of $\L(H)$, with 
$X\subset Y$; direct $\cS$ by inclusion. 
For each $\alpha\in\cS$, let $P_\alpha :Y\to X$ be a complete projection of 
$Y$ onto $X$, with $\|P_\alpha\|_{cb} \le\lambda$. 
We may now use the reflexivity of $Y$ and the Tychonoff theorem to produce 
a completely bounded projection from $\L(H)$ onto $Y$. 
For each $\alpha\in \cS$, define $\tilde P_\alpha : \L(H) \to X$ by 
$\tilde P_\alpha (v) = 0$ if $v\notin \alpha$; $\tilde P_\alpha (v)=
P_\alpha (v)$ if $v\in \alpha$. 
$\tilde P_\alpha$ is neither continuous nor linear; nevertheless, the 
weak-compactness of the ball of $X$ yields a subnet 
$(P_{\alpha_\beta})_{\beta\in \D}$ of $(P_\alpha)_{\alpha\in\cS}$ so that 
$Pv \defeq \lim_{\beta\in\D} P_{\alpha_\beta}(v)$ exists weakly for all 
$v\in B(H)$. 
Since of course every $v\in B(H)$ is contained in some $\alpha\in\cS$, it 
follows that $P$ is indeed a linear projection from $\L(H)$ into $X$. 
Finally, we also have that for all $n$, 
$$\|P\|_n \le \varlimsup_{\alpha\in\cS} \|P_\alpha\|_n\le \lambda\ ,$$ 
whence $\|P\|_{cb} \le\lambda$, showing that $X$ is indeed 
$\lambda$-injective.\qed  
\enddemo

\remark{Comment} 
This (rather outrageous) use of the Tychonoff theorem is due to 
J.~Lindenstrauss.
\endremark 

\remark{Remark} 
We show in Proposition 2.22 below that if $X$ is a separable operator space 
with the {\rm $\lambda$-CSEP} and $\lambda<2$, then $X$ is reflexive 
and (hence is $\lambda$-injective).
\endremark 

Work of G.~Pisier's yields immediately that every separable reflexive 
operator space which is isomorphically injective is Hilbertian, i.e., 
Banach isomorphic to Hilbert space (cf. \cite{R}). 
Evidently Corollary~2.9 also yields that every isomorphically injective 
separable operator space is reflexive (and so Hilbertian). 

Of course the natural (and far from obvious!) special problem in this setting 
is then as follows: 
\roster
\item"{}" {\it let $X$ be a separable infinite-dimensional isomorphically 
injective operator space. 
Is $X$ completely isomorphic to $\bR$, $\bC$, or $\bR\oplus \bC$?} 
\endroster 
(This problem has been solved affirmatively for 
$X$ isometrically injective by A.~Robertson \cite{R}.)  
A remarkable result of T.~Oikhberg \cite{O} yields the answer is 
{\it affirmative\/} if $X$ is completely isomorphic to a subspace of 
$\bR\oplus \bC$. 
Finally, we note the following quantitative problem, 
whose positive solution implies an affirmative answer to 
the preceding question, in virtue of Oikhberg's result.  
\roster
\item"{}" {\it Let $X$ be a separable operator space which is 
$\lambda$-injective. 
Is there a $\beta$, depending only on $\lambda$, and an $n$ (depending 
on $X$) so that $d_{cb}(X,Y)\le\beta$ for some subspace $Y$ of $M_{\infty,n}
\oplus M_{n,\infty}$?} 
\endroster

Before dealing with the main result of this section, we give the 

\demo{Proof of Proposition 2.6} 
We identify $M_{\infty,n}$ with $\bC\otimes \bR_n$ endowed with its natural 
operator space structure (where $R_n$ denotes the $n$-dimensional row space). 
Letting $e_1,\ldots,e_n$ be the natural orthonormal basis of $\bR_n$, any 
vector $v\in \bC\otimes \bR_n$ has the form 
$$v = \sum_{i=1}^n u_i \otimes e_i\ \text{ for unique }\ u_1,\ldots,u_n\in 
\bC\ .
\tag 2.2$$
In fact the map $P_i$ which sends $v$ to $u_i$, yields a projection from  
$\bC\otimes \bR_n$ onto $\bC$. 
Now letting $F$ be a $k$-dimensional subspace of $\bC\otimes \bR_n$ and 
setting $V_i = PF_i$ we have that $V_i$ is a subspace of $\bC$ with 
$\dim V_i\le k$ for all $i$, and clearly 
$$F\subset \mathop{\text{span}}_{1\le i\le n} V_i\otimes e_i \subset 
\bV \otimes \bR_n 
\tag 2.3$$
where $\bV = V_1+\cdots + V_n$. 
Evidently $m\defeq \dim \bV \le k\cdot n$; by homogeneity of $\bC$, 
$\bV \otimes \bR_n$ is completely isometric to $M_{m,n}$, which in turn is 
completely isometric to a subspace of $M_{k\cdot n}$. 
The proof for $M_{n,\infty}$ is of course the same.\qed
\enddemo 

We now prove Theorem 2.2, giving the full details of the modification of our 
original argument, due to T.~Oikhberg. 
We shall see the argument is essentially the same as the one alluded to in 
Section one, after inserting the appropriate quantizations. 
We first give the reformulation  analogous to Theorem~1.5. 

\proclaim{Theorem 2.11} 
Let $Z_1,Z_2,\ldots$, $X$, and $Y$ be operator spaces satisfying the hypotheses 
of Theorem~2.2. 
Let $(T_j)$ be a sequence of completely bounded operators, with $T_j:Y\to Z_j$ 
for all $j$, so that $T_j|X\to 0$ in the {\rm SOT}, and $\sup_j\|T_j\|_{cb}=1$. 
For every $\ep >0$, there exists a sequence $(S_j)$ of completely bounded 
operators so that for all $j$, 
\roster
\item"1)" $S_j :Y\to Z_j$ 
\item"2)" $X\subset \ker S_j$ 
\item"3)" $\|S_j\|_{cb} < C \lambda + {\ep\over2\lambda}$ 
\endroster 
so that $(T_j-S_j)\to 0$ in the {\rm SOT}.
\endproclaim 

We first give the proof that Theorem~2.11 $\Rightarrow$ Theorem~2.2. 
Let $T:X\to Z$ be as in the statement of 2.2, and let $(\tilde T_j)$  be the 
sequence so that $TX = (\tilde T_jx)$ for all $x\in X$. 
For each $j$, since $Z_j$ is $\lambda$-injective, we may chose 
$T'_j :Y\to Z_j$ with $\|T'_j\|_{cb} \le \lambda\|\tilde T_j\|_{cb}$ and 
$T'_j\mid X= \tilde T_j$. 
Of course then 
$$\sup_j \|T'_j\|_{cb} \defeq \beta \le \lambda\|T\|_{cb}\ .
\tag 2.4$$ 
Let $T_j = T'_j/\beta$ for all $j$. 
Now the hypotheses of 2.11 hold; choose $(S_j)$ satisfying its conclusion. 
Now define $\ttT :Y\to (Z_j)_{\ell^\infty}$ by 
$$\ttT (y) = (T_j-S_j)(y)\ \text{ for all }\ y\in Y\ . 
\tag 2.5$$

It follows from the conclusion of 2.11 that $\ttT$ actually has its range in 
$Z= (Z_j)_{c_0}$, and 2) of 2.11 insures that $\ttT$ extends $\frac1{\beta}T$. 
Thus defining $\tilde T = \beta \ttT$, $\tilde T$ extends $T$ and by (4) 
and 3) of 2.11 
$$\align 
\|\tilde T\|_{cb} & < \beta \left(1+C\lambda +{\ep\over2}\right)\tag 2.6\cr
&< (C\lambda^2 + \lambda +\ep)\|T\|_{cb}\ .
\endalign$$
\vskip-30pt
\rightline{$\square$}
\medskip

We again need the analogue for Lemma~1.6, which actually holds with no 
uniform exactness assumption. 

\proclaim{Lemma 2.12} 
Let $Z_1,Z_2,\ldots$ be arbitrary operator spaces, $X\subset Y$ operator 
spaces  with $Y/X$ finite-dimensional, and $(T_j)$ a sequence of completely 
bounded operators so that for all $j$, $T_j :Y\to Z_j$ with 
$\|T_j\|_{cb}\le1$, so that $T_j|X\to 0$ in the {\rm SOT}. 
Let $F$ be a finite-dimensional subspace of $Y$ with $X\oplus F=Y$, and 
let $P$ be the projection of $Y$ onto $F$ with kernel $X$. 
Then for all positive integers $n$, 
$$\varlimsup_{j\to\infty} \|T_j P\|_n\le 1\ .
\tag 2.7$$
\endproclaim 

\demo{Proof} 
Suppose not. 
Then by passing to a subsequence if necessary, we can fix an $n$ and wlg 
choose $C>1$ and $(A_j)$ a norm-one sequence in 
$M_n(Y) = \L(\ell_n^2)\otimes Y$ with 
$$\|(I\otimes T_j) \cdot (I\otimes P)(A_j) \| > C\ \text{ for all }\ j\ . 
\tag 2.8$$ 
Now setting 
$$A_j = \left[ \matrix 
y_{11}^j & \cdots & y_{1n}^j\cr 
\vdots&&\vdots\cr 
y_{n1}^j&\cdots &y_{nn}^j 
\endmatrix \right]
\tag 2.9$$ 
for each $j$, choose $U_j \in M_n(X)$ and $V_j\in M_n(F)$ so that 
$$A_j = U_j \oplus V_j\ ,
\tag 2.10$$ 
whence by (8) and the fact that $(I\otimes P)(A_j) = V_j$; 
$$\|(I\otimes T_j)(V_j)\| > C\ . 
\tag 2.11$$ 
Of course for each $j$, we may set 
$$U_j = \left[ \matrix x_{11}^j&\cdots & x_{1n}^j\cr 
\vdots&&\vdots\cr 
x_{n1}^j&\cdots&x_{nn}^j \endmatrix \right] 
\tag 2.12i$$
and 
$$V_j = \left[ \matrix f_{11}^j&\cdots & f_{1n}^j\cr 
\vdots&&\vdots\cr 
f_{n1}^j&\cdots&f_{nn}^j \endmatrix \right] \ .
\tag 2.12ii$$

Thus (11) means that 
$$\left\| \left[ \matrix T_jf_{11}^j &\cdots & T_j f_{1n}^j\cr 
\vdots&&\vdots\cr
T_jf_n^j&\cdots & T_j f_{nn}^j\endmatrix \right] \right\| 
> C\ \text{ for all } j\ .
\tag 2.13$$ 
Since $F$ is a finite-dimensional space, $P$ is bounded and of course 
completely bounded; in particular, the sequences $(U_j)$ and $(V_j)$ are 
both bounded, so by compactness of bounded subsets of $M_n(F)$, by passing 
to a further subsequence if necessary, we may assume for some $V\in M_n(F)$ 
that 
$$V_j \to V\ \text{ in norm.} 
\tag 2.14$$
(In other words, we have 
$$V= \left[ \matrix f_{11}&\cdots&f_{1n}\cr \vdots&&\vdots\cr 
f_{n1}&\cdots&f_{nn}\endmatrix\right] $$
and for each $i$ and $k$, $f_{ik}^j \to f_{ik}$ in norm.) 
But then 
$$\|(I\otimes T_j) (V_j-V)\|\to 0 
\tag 2.15$$ 
whence by (2.11) 
$$\varliminf_{j\to\infty} \|(I\otimes T_j)(V)\| \ge C\ .
\tag 2.16$$
That is 
$$\varliminf_{j\to\infty} \left\| \matrix 
T_j f_{11}&\cdots&T_j f_{1n}\cr
\vdots&&\vdots\cr
T_jf_{n1}&\cdots&T_j f_{nn}\endmatrix \right\| \ge C \ .$$
Now also 
$$\|U_j\oplus V\| \to 1
\tag 2.17$$
since $\|U_j \oplus V_j\| = \|A_j\| =1$ for all $j$ and by (2.14), 
$\|U_j \oplus V_j\| - \|U_j \oplus V\| \to0$. 
Now fix $\ep >0$ with $1+\ep <C$, and choose $k$ with 
$$\|U_k\oplus V\| < 1+\ep\ \text{ (using (2.17)).}
\tag 2.18$$

Then for all $j$, 
$$\|(I\otimes T_j)(U_k\oplus V)\| \le \| T_j\|_{cb} \| U_k\oplus V\| 
< 1+\ep\ .
\tag 2.19$$
Since $T_j |X \to 0$ in the SOT, $I\otimes T_j|M_n(X)\to 0$ in the SOT, whence 
$$\lim_{j\to\infty} \|(I\otimes T_j)(U_k\oplus V)\| 
- \|I\otimes T_j)(V)\| =0 
\tag 2.20$$
Thus we obtain from (19) and (20) that 
$$\varlimsup_{j\to\infty}\|(I\otimes T_j)(V)\| \le 1+\ep \ ,
\tag 2.21$$ 
contradicting (2.16).\qed
\enddemo 

We now  apply a useful result of R.~Smith to obtain the following consequence. 

\proclaim{Corollary 2.13} 
Let $Z_1,Z_2,\ldots,$  $X$, $Y$, and $(T_j)$, $F$ and $P$ be as in Lemma~2.12, 
and assume $\{Z_1,Z_2,\ldots\}$ is $C$-uniformly exact. 
Then 
$$\varlimsup_{j\to\infty} \|T_j P\|_{cb} \le C\ .
\tag 2.22$$ 
\endproclaim 

\demo{Proof} 
Roger Smith's lemma \cite{S} yields that for all $n$, operator spaces $X$, and 
linear maps $T:X\to M_n$, 
$$\|T\|_{cb} = \|T\|_n\ .
\tag 2.23$$ 
(See \cite{Pi} for the operator space formulation and another proof.) 
Let $k=\dim F$ and $n= \bn (k)$ where $\bn$ is the $C$-uniformity function 
for $\{Z_1,Z_2,\ldots\}$. 
Then fixing $j$, the range of $T_jP$ is a subspace of $Z_j$ of 
dimensional at most $n$, hence we obtain from (2.23) and Definition~2.3 
that 
$$\|T_jP\|_{cb} \le C\|T\|_n\ ,
\tag 2.24$$ 
which immediately yields (2.22) in virtue of Lemma~2.12.\qed
\enddemo 

The next result follows from this Corollary in the same manner as the 
corresponding 
Banach space result follows from Lemma~1.6 (see the Remark following the 
proof of 1.6). 

\proclaim{Corollary 2.14} 
Let $X\subset Y$ be operator spaces with $Y/X$ finite-dimensional, $(Z_j)$ 
a sequence of operator spaces so that $\{Z_1,Z_2,\ldots\}$ is $C$-uniformly 
exact, and $T:Y\to (Z_j)_{\ell^\infty}$ a non-zero completely bounded 
operator with $TX \subset (Z_j)_{c_0}$. 
Then given $\ep>0$, there exists $\tilde T:Y\to (Z_j)_{c_0}$ so that 
$\tilde T$ extends $T|X$ and $\|\tilde T\|_{cb}\le (C+1+\ep)\|T\|_{cb}$. 
\endproclaim 

\demo{Proof of Theorem 2.11} 
Assume (in virtue of Corollary 2.14) that $Y/X$ is infinite-dimensional, and 
let $y_1,y_2,\ldots$, $F_k$, $Y_k$ and $P_k$ be as in the proof of 
Theorem~1.5; let $\ep>0$. 
We construct for each $k$, a sequence $(S_n^{(k)})$ of operators 
so that there is an $M_k$ such that for all $n$, 
$$\gather
\text{if } n<M_k, \text{ then } S_n^{(k)} :Y\to Z_n \text{ and } 
\|S_n^{(k)}\|_{cb} < C\lambda + {\ep\over2\lambda} 
\tag 2.25\\
\text{if } n\ge M_k,\ S_n^{(k)} :Y_k \to Z_n,\ S_n^{(k)} = T_nP_k\text{ and } 
\|S_n^{(k)}\|_{cb} < C+{\ep\over 2\lambda^2} 
\tag 2.26\\ 
X\subset \ker S_n^{(k)} \tag 2.27\\
S_n^{(k+1)} | Y_k = S_n^{(k)} \ .\tag 2.28
\endgather $$
Letting $S_n = \overline{\bigcup_{k=1} S_n^{(k)}}$, it follows that $S_n$ 
is a well-defined completely bounded operator with $\|S_n\|_{cb} \le 
C\lambda + {\ep\over2\lambda} < C\lambda +{\ep\over\lambda}$ for all $n$. 
Just as before, it then follows that $(S_n)$ satisfies the conclusion of 
Theorem~2.11. 

Again, we construct the sequences $(S_n^{(k)})_{n=1}^\infty$ by induction  
on $k$, setting $S_n^{(0)} = 0$ for all $n$. 
Let $k\ge0$ and suppose $S_n^{(k)}$, $M_k$ have been defined, satisfing 
(2.25)--(2.27). 
By Corollary~2.13, choose $M_{k+1} > M_k$ so that for all $n\ge M_{k+1}$, 
$$\|T_nP_{k+1}\|_{cb} < C+ {\ep\over 2\lambda^2}\ .
\tag 2.29$$ 
Now for $n\ge M_{k+1}$, let $S_n^{(k+1)} = T_n P_{k+1}$. 
Since also $S_n^{(k)} = T_n P_k$ by (2.26), we have that (2.28) holds. 

Now for $M_k\le n< M_{k+1}$, since $Z_n$ is a $\lambda$-injective operator 
space, choose $S_n^{(k+1)}$ a linear extension of $S_n^{(k)}$ from 
$Y_k$ to $Y$ with 
$$\|S_n^{(k+1)}\|_{cb} \le \lambda \|S_n^{(k)}\|_{cb} < C\lambda 
+ {\ep\over 2\lambda}\ \text{ (by (2.26))}.
\tag 2.30$$ 
Finally, for $n<M_k$, let $S_n^{(k+1)} = S_n^{(k)}$. 
Evidently (2.28) holds for all $n< M_{k+1}$, and of course we have that 
(2.25)--(2.27) hold replacing ``$k$'' by ``$k+1$''. 
(Again the procedure is also valid in the setting of Corollary~2.5, since 
the separability of $Y$ and the assumption that $Z_n$ has the $\lambda$-CSEP 
allows us to do this.)\qed  
\enddemo 

Theorem 2.11 (and hence Theorem 2.2) holds under a hypotheses weaker than that 
of uniform exactness. 
Here is the relevant concept. 

\definition{Definition 2.4} 
A family $\Z$ of operator spaces is said to be of {\it finite matrix type\/} 
if there is a $C\ge1$ so that for any finite-dimensional operator space $F$, 
there is an $n= \bn (F)$ so that 
$$\|T\|_{cb} \le C\|T\|_n\text{ for all linear operators } 
T:F\to Z\text{ and all } Z\in \Z\ .
\tag 2.31$$ 
If $C$ works, we say that $\Z$ is of $C$-finite matrix type, or briefly, 
$\Z$ is $C$-finite; if the function $\bn$ works, we say that $\Z$ is 
$C$-finite with function $\bn$. 
\enddefinition

\noindent 
(Note that the domain of $\bn$ is the family of all finite-dimensional 
operator spaces.) 
An operator space $Z$ is $C$-finite provided $\{Z\}$ is $C$-finite. 

Thanks to the result of R.~Smith cited in the proof of Corollary~2.13, it 
follows that if $\Z$ is $C$-uniformly exact, $\Z$ is $C$-finite. 

\proclaim{Proposition 2.15} 
Theorems 2.2 and 2.11 both hold if one replaces the assumption that 
$\{Z_1,Z_2,\ldots\}$ is $C$-uniformly exact by the assumption that 
$\{Z_1,Z_2,\ldots\}$ is $C$-finite. 
\endproclaim 

Of course it follows that also Corollaries 2.4 and 2.5 hold under this 
weaker assumption. 
T.~Oikhberg has actually obtained a converse to this result, which goes 
as follows (see \cite{OR}): 
{\it Let $Z_1,Z_2,\ldots$ be separable operator spaces so that $(Z_j)_{c_0}$ 
has the {\rm CSEP}. 
Then $\{Z_1,Z_2,\ldots\}$ is of finite matrix type\/}. 

Proposition 2.15 follows readily from the argument for 2.11 and the following 
simple, useful result about operator spaces: 

\proclaim{Fact} 
Let $Y,Z$ be operator spaces and $T:Y\to Z$ be a completely bounded map; 
let $X$ be a (closed linear) subspace of $Y$, 
let $\ker T\supset X$, $\pi:Y\to Y/X$ the quotient map, and 
$\tilde T :Y/X\to Z$ the canonical map with $T= \tilde T\pi$. 
Then  
\roster 
\item"(a)" $\|T\|_n = \|\tilde T\|_n$ for all $n$
\endroster
and hence
\roster
\item"(b)" $\|T\|_{cb} = \|\tilde T\|_{cb}$. 
\endroster
\endproclaim

\noindent 
(This fact in turn follows from the natural, elementary result: 
if $n\ge1$ and $(y_{ij})\in M_n(Y)$ (regarded as contained in $M_{00} (Y)$), 
then $\|\pi (y_{ij})\| = \inf \|(y_{ij})-(x_{ij})\| :(x_{ij}) \in M_n(X)$). 

\demo{Proof of Proposition 2.15} 
The main new observation is that the conclusion of Corollary~2.13 holds 
if we assume instead that $\{Z_1,Z_2,\ldots\}$ is $C$-finite. 
Indeed, let $G= Y/X$ and let $n=\bn (G)$, where $\bn$ is as in 
Definition~2.4. 
Now fix $j$, let $\pi :Y\to Y/X$ be the quotient map, and 
$\tilde T_j:G\to Z_j$ the canonical map with $\tilde T_j\pi= T_jP$. 
Then 
$$\align 
\|T_jP\|_{cb} & = \|\tilde T_j\|_{cb}\ \text{ by the Fact}\cr
&\le C\|\tilde T_j\|_n\ \text{ since $Z_j$ is $C$-exact}\cr 
&= C\|T_jP\|_n\ \text{ by the Fact.}
\endalign$$ 
Hence we deduce from Lemma 2.12 that 
$$\varlimsup_{j\to\infty} \|T_jP\|_{cb} \le C$$ 
as desired. 
The proof of the modified statement of Theorem~2.11 is now identical to the 
argument given above, whence 2.15 follows.\qed
\enddemo 

We finally deal with certain  quantized formulations of the later 
results of Section~1. 
We first note that Corollary~2.4 follows from Theorem~2.11 in exactly the 
same way as Corollary~1.7 follows from Theorem~1.5. 
We also have the following quantized form of local complementability 
(Definition~1.3). 

\definition{Definition 2.5} 
Let $X\subset Y$ be given operator spaces. 
$X$ is said to be {\it completely locally complemented\/} in $Y$ if there is a 
$\lambda \ge 1$ so that 
$$ \eqalign{
&X\text{ is a completely $\lambda$-complemented in $Z$ for all}\cr
&X\subset Z\subset Y\text{ with } Z/X\text{ finite-dimensional.}\cr} 
\tag 2.32$$ 
When (2.32) holds, we say $X$ is completely locally $\lambda$-complemented 
in $Z$. 
We also say that $X$ is completely locally $\lambda$-cocomplemented in $Y$ 
if for all $X\subset Z\subset Y$ with $Z/X$ finite-dimensional there exists 
a projection $P:Z\to X$ with $\|I-P\|_{cb} \le\lambda$.
\enddefinition 

Proposition 2.15, Corollary 2.14 and the arguments for Corollary 1.8 now 
immediately yield the following result. 

\proclaim{Corollary 2.16} 
Let $Z_1,Z_2,\ldots$ be a sequence of operator spaces so that $\{Z_1,Z_2,
\ldots\}$ is $C$-finite, and let $Z= (Z_j)_{c_0}$. 
Then $Z$ is completely locally $(C+\ep)$-cocomplemented in 
$(Z_j)_{\ell^\infty}$. 
\endproclaim 

The Corollary in the Remark following Lemma~1.6 also immediately yields 

\proclaim{Corollary 2.17} 
Let $Z_1,Z_2,\ldots$ be a sequence of one-injective operator spaces and 
again let $Z= (Z_j)_{c_0}$. 
Then $Z$ is Banach $(2+\ep)$-locally complemented in any operator superspace. 
\endproclaim

\remark{Remark} 
Corollaries 2.16 and 2.17 have also been obtained in a different way in 
\cite{JO}. 
\endremark 

The concept of local complementability may be used to refine the formulation 
of Proposition~2.10, as follows. 

\proclaim{Proposition 2.18} 
Let $X$ be a reflexive operator space, and suppose that $X$ is completely 
locally $\lambda$-complemented in every operator space $Y$ with $X\subset Y$. 
Then $X$ is $\lambda$-injective.
\endproclaim 

\demo{Proof} 
The proof is really the same as  the argument for 2.10. 
We assume that $X\subset \L(H)$ for $H$ a (not necessarily separable) Hilbert 
space, and then let $\cS$ be the family of subspaces $Y$ of $H$ with 
$X\subset Y$ and $\dim Y/X<\infty$. 
$\cS$ is again directed by inclusion, and the argument that $X$ is 
completely $\lambda$-complemented in $\L(H)$ now follows just as before.\qed
\enddemo 

We finally note the quantized versions of Theorem~1.9 and its consequences. 

\proclaim{Theorem 2.19} 
Replace ``bounded'' by ``completely bounded'' and ``$\|\tilde T\| \le 
2\|T\|$'' by $\|\tilde T\|_{cb} \le 2\|T\|_{cb}$ in the statement of 
Theorem~1.9. 
\endproclaim 

\demo{Proof} 
Assume wlg that $\|T\|_{cb} =1$. 
Let $\M$ be as in Lemma~1.10, and define $\cS$ by 
$$\cS = \{S\in \L(Y,Z) \mid \|S\|_{cb}\le 1\text{ and } X \subset\ker S\}. 
\tag 2.33$$ 
Again, we see that all SOT-cluster points of $(T_n)$ lie in $\cS$, for if 
$T_{n_i} \to S$ in the SOT, also for each fixed $k$, 
$\|T_{n_i}-S\|_k\to 0$, whence $\|T_{n_i}\|_k \to \|S\|_k$, so 
$\|S\|_k\le 1$, and thus $\|S\|_{cb} \le 1$. 
Of course also $S(X)=0$ for all $x\in X$; so $S\in \cS$. 
This argument also shows that $\cS$ is SOT-closed. 
The rest of the argument is the same as for 1.9.\qed
\enddemo 

The next result follows again immediately from 2.19 and the proof of 
Corollary~1.12.

\proclaim{Corollary 2.20} 
Let $Z$ be a finite-dimensional operator space and $Y$ a separable subspace 
of $\ell^\infty (Z)$ containing $c_0(Z)$. 
Then $c_0(Z)$ is completely contractively cocomplemented in $Y$.
\endproclaim 

\proclaim{Corollary 2.21} 
$c_0(M_n)$ has the {\rm 2-CSEP}, for all $n$.
\endproclaim 

This follows immediately from 2.20 and the fact that $M_n$ is a 
1-injective operator space. 
The next section gives a ``saving property'' for the space $(M_n)_{c_0}$ 
in view of its failure to have the CSEP. 

Our last result of this section shows that separable non-injective operator 
spaces cannot have the $\lambda$-CSEP if $\lambda<2$. 

\proclaim{Proposition 2.22} 
Let $X$ be a separable operator space with the {\rm $\lambda$-CSEP}. 
If $\lambda<2$, then $X$ is reflexive (and hence is $\lambda$-injective 
by Proposition~2.10).
\endproclaim  

We first require the corresponding ``pure'' Banach space result. 
(This  may be part of the subject's folk-lore. 
The argument I give is due to W.B.~Johnson, and I'm most grateful to him 
for providing this elegant proof.) 

\proclaim{Lemma 2.23} 
Let $X$ be a separable Banach space containing a subspace $Y$ isomorphic 
to $c_0$. Given $\ep>0$, there exists a subspace $Z$ of $Y$ so that $Z$ 
is $(1+\ep)$-isomorphic to $c_0$ and $Z$ is $(1+\ep)$-complemented in $X$.
\endproclaim 

\demo{Proof} 
Let $\bep >0$ be such that $(1+\bep)^2 < 1+\ep$.
By a result of R.C.~James \cite{J}, we may chose a subspace $E$ of $Y$ with 
$E$ $(1+\bep)$-isomorphic to $c_0$. 
It follows that we may choose a basis $(e_j)$ for $E$ so that for all 
$(c_j)$ in $c_0$, 
$$\sup_j \sum |c_j| \le \Big\|\sum c_j e_j\Big\| \le (1+\bep) \sup |c_j|\ .
\tag 2.34$$ 
Let $(f_n)$ be a Hahn-Banach extension to $X$ of the biorthogonal functionals 
to $(e_n)$. 
By passing to a subsequence, we may assume w.l.g. that $(f_n)$ converges 
$w^*$ in $X^*$. 
Now define $g_n$ by 
$$g_n = {f_{2n} - f_{2n-1} \over 2} \ .
\tag 2.35$$ 
It follows that $g_n\to0$ $w^*$ and moreover (by (2.34)), 
$$\|g_n\| < 1+\bep \ \text{ for all }\ n\ . 
\tag 2.36$$ 
Finally, let $z_n = e_{2n} - e_{2n-1}$ for all $n$; then let $Z= [z_n]$. 
Of course $Z$ is $(1+\ep)$-isomorphic to $c_0$, and $(g_n)$ is biorthogonal 
to $(z_n)$. 
Thus we may define a projection $P:X\to Z$ by 
$$Px = \sum g_n(x) z_n\ \text{ for all } x\in X\ . 
\tag 2.37$$ 
It follows that if $x\in X$, then 
$$\align 
\|Px\| & \le (1+\bep) \sup_n | g_n(x)| \  \text{ by (2.34)}\tag 2.38\cr 
&\le (1+\bep)^2 \|x\|\qquad\quad \text{ by (2.36)}\ .
\endalign$$
Hence $P$ is indeed a projection onto $Z$ with $\|P\| < 1+\ep$.\qed
\enddemo

\demo{Proof of Proposition 2.22}
Suppose to the contrary that $X$ is not reflexive. 
Then $X$ contains a subspace isomorphic to $c_0$ by Corollary~2.9. 
Now let $\ep>0$, to be decided later, and choose by Lemma~2.23 a subspace $Z$ 
of $X$ which is (Banach) $(1+\ep)$-isomorphic to $c_0$ and 
$(1+\ep)$-complemented in $X$. 
Now let $Y$ be a separable subspace of $Z^{**}$ with $Z\subset Y$, let 
$i:Z\to X$ be the identity injection, and also let $P:X\to Z$ be a projection 
with $\|P\| <1+\ep$. 
Since $X$ has the $\lambda$-CSEP, letting $Y$ have its natural operator 
space structure, we find a completely bounded extension $\tilde \imath :
Y\to X$ with $\|\tilde\imath\|_{cb} \le\lambda$. 
But then letting $Q= P\tilde \imath$, $Q$ is a projection from $Y$ onto 
$Z$ and 
$$\|Q\| < (1+\ep)\lambda\ .
\tag 2.39$$ 
Since $Z$ is $(1+\ep)$-isomorphic to $c_0$, it now follows that if 
$\tilde Y$ is separable with $c_0 \subset \tilde Y\subset \ell^\infty$, then 
$$c_0\ \text{ is } (1+\ep)^2 \lambda\text{-complemented in }\tilde Y\ .
\tag 2.40$$
But this implies that $c_0$ itself has the $(1+\ep)^2\lambda$-SEP, hence 
by Sobczyk's result \cite{S}, $(1+\ep)^2\lambda \ge2$. 
Of course we then need only choose $\ep>0$ with $(1+\ep)^2\lambda <2$, to 
arrive at the desired contradiction.\qed
\enddemo 

\head 3. The Complete  Separable Complementation Property\endhead 

In this section we study the following concept, more general than the CSEP.

\definition{Definition 3.1} 
A separable operator space $Z$ has {\it Complete Separable Complementation 
Property\/} (the CSCP) if whenever $Y$ is a separable 
locally reflexive operator space, $X$ is a subspace of $Y$, and $T:X\to Z$ 
is a complete surjective isomorphism, $T$ has a completely bounded 
extension $\tilde T:Y\to Z$.
\enddefinition 

In other words, $Z$ has the CSCP 
provided every complete isomorph of 
$Z$ is completely complemented in every separable locally reflexive 
operator superspace. 

\remark{Remark} 
After the first draft of this paper was completed, it was discovered that 
if $Z$ has the CSCP, then the diagram (0.1) holds for {\it arbitrary\/} 
completely bounded maps $T$; moreover $Z$ has the CSCP provided it is 
completely complemented in every locally reflexive separable operator 
superspace (see \cite{OR}). 
\endremark 

Evidently this property is invariant under complete isomorphisms. 
The main result of this section is as follows. 

\proclaim{Theorem 3.1} 
Let $\lambda\ge1$ and 
let $Z_1,Z_2,\ldots$ be separable $\lambda$-injective operator spaces. 
Then $(Z_j)_{c_0}$ has the {\rm CSCP}. 
\endproclaim 

\proclaim{Corollary} 
$\bK_0$ has the {\rm CSCP}.
\endproclaim 

In fact, our proof yields that 
{\it if $Y$ is a $C$-locally reflexive separable superspace of $(Z_j)_{c_0}$, 
then $(Z_j)_{c_0}$ is completely 
$\lambda^3+ (C+1)\lambda^2+\lambda+\ep$-complemented in $Y$, for all 
$\ep>0$.} 
As in the preceding section, this Theorem follows via the modification 
of T.~Oikhberg of our original construction for the case $\lambda=1$. 
Because of the known structure of the separable isometric injectives 
(\cite{R}), 
Theorem~3.1 for the case $\lambda =1$ is equivalent to: 
$(M_{\infty,n}\oplus M_{n,\infty})_{c_0}$ {\it has the\/} CSCP. 
After the first draft of this paper was completed, it was discovered that 
$\bK$ (the space of 
compact operators on $\ell^2$), {\it has\/} the CSCP. 
(The proof uses the above Corollary --- see \cite{OR}.) 
The main structural problem for this property is as follows: 

\proclaim{Problem}
Is every space with the {\rm CSCP} completely isomorphic to a subspace 
of $\bK$?
\endproclaim 

We discuss further aspects of this problem in Section~4. 
Let us also note that by T.~Oikhberg's result (see \cite{OR}), 
Theorem~3.1 {\it fails\/} without the assumption  of local reflexivity 
in the definition of the CSCP. 
Positive motivation for Theorem~3.1 and Definition~3.1 is given by 
the following result: 

\proclaim{Corollary 3.2} 
Let $Z_1,Z_2,\ldots$ be as in the statement of Theorem~3.1, $\A$ be a 
separable nuclear $C^*$-algebra, and $\tilde Z$ be a subspace of $\A$ 
which is completely isomorphic to $(Z_j)_{c_0}$. 
Then $\tilde Z$ is completely complemented in $\A$.
\endproclaim 

This follows immediately from Theorem 3.1, in virtue of the fact that 
nuclear $C^*$-algebras are 1-locally reflexive \cite{EH}. 
The quantitative version of 3.1 yields 

\proclaim{Corollary 3.3} 
Let $\A$ be a separable nuclear $C^*$-algebra and $\K_0$ be a 
$C^*$-algebra which is *-isomorphic to $\bK_0 = (M_n)_{c_0}$. 
Then for all $\ep>0$, $\K_0$ is completely $(4+\ep)$-complemented in $\A$.
\endproclaim 

Corollary 3.2 also suggests the following 

\proclaim{Problem} 
Let $Z$ be a separable operator space which completely embeds in some nuclear 
$C^*$-algebra. 
Suppose that every complete embedding of $Z$ into a nuclear separable 
$C^*$-algebra $\A$ is completely complemented in $\A$. 
Does $Z$ have the {\rm CSCP}?
\endproclaim 

We now deal with the proof of Theorem 3.1. 
As has been the case in the preceding section, the arguments hold in 
considerable generality; {\it local complementability\/} (cf.\ Definition~2.5) 
plays a key role in the discussion. 

\proclaim{Theorem 3.4} 
Let $\lambda\ge1$, and 
let $Z_1,Z_2,\ldots$ be $\lambda$-injective operator spaces and $X\subset Y$ 
operator spaces with $X$ locally complemented in $Y$ and $Y/X$ separable. 
Let $T:X\to Z$ be a completely bounded operator, where $Z= (Z_j)_{c_0}$. 
Then $T$ admits a completely bounded extension $\tilde T:Y\to Z$. 
\endproclaim

\remark{Remark} 
The proof yields that if $X$ is completely locally 
$\beta$-cocomplemented in $Y$, 
then for all $\ep>0$, $\tilde T$ may be chosen with $\|\tilde T\|_{cb} 
< (\beta \lambda^2 +\lambda +\ep)\|T\|_{cb}$ (if $T\ne0$). 
\endremark

We first note an immediate consequence. 

\proclaim{Corollary 3.5} 
Let $Z_1,Z_2,\ldots$ and $Z$ be as in Theorem~3.4 and let $Y$ be an 
operator space with $Z\subset Y$, $Y/Z$ separable, and $Z$ locally 
complemented in $Y$. 
Then $Z$ is completely complemented in $Y$. 
\endproclaim 

\remark{Remark} 
Again, we obtain that if $Z$ is completely locally $\beta$-cocomplemented 
in $Y$, $Z$ is completely 
$(\beta\lambda^2 +\lambda +\ep)$-complemented in $Y$ for all 
$\ep>0$.
\endremark 

As before, we first reformulate Theorem 3.4. 

\proclaim{Theorem 3.6} 
Let $Z_1,Z_2,\ldots$, $Z$, $X$, and $Y$ be operator spaces satisfying 
the hypotheses of Theorem~3.4, and suppose $X$ is completely locally 
$C$-cocomplemented in $Y$. 
Let $(T_j)$ be a sequence of completely bounded operators with $T_j:Y\to Z_j$ 
for all $j$, so that $T_j|X\to0$ in the {\rm SOT} and $\sup_j\|T_j\|_{cb}=1$. 
{For} every $\ep>0$, there exists a sequence $(S_j)$ of completely bounded 
operators so that for all $j$, 
\roster
\item"1)" $S_j:X\to Z_j$ 
\item"2)" $X\subset \ker S_j$
\item"3)" $\|S_j\|_{cb} < C\lambda + {\ep\over 2\lambda}$
\endroster
so that $(T_j-S_j)\to0$ in the {\rm SOT}. 
\endproclaim 

The proof that Theorem 3.6 $\Rightarrow$ Theorem 3.4 is again the same 
as the one showing Theorem~2.11 $\Rightarrow$ Theorem~2.2; this proof also 
yields the quantitative statement in the Remarks following Theorem~3.6, 
as well as the following quantitative variation
of Corollary~3.5 (all objects as 
in its statement): 
{\it If $Z$ is completely locally $C$-cocomplemented in $Y$, then $Z$ 
is $(C\lambda +\ep)$-cocomplemented in $Y$ for all $\ep>0$.} 

The proof of Theorem 3.4 (i.e., of Theorem 3.6) is analogous to the proofs  
of Theorems~1.1 and 2.2; it requires a different (again rather surprising) 
lemma, replacing Lemmas~1.6 and 2.12. 

\proclaim{Lemma 3.7} 
Let $Z_1,Z_2,\ldots$, $X$, $Y$ be arbitrary operator spaces with $X\subset Y$ 
and $Y/X$ finite-dimensional. 
Let $(T_n)$ be a sequence of completely bounded operators with $T_j:Y\to Z_j$ 
for all $j$, so that  $T_n|X\to0$ in the {\rm SOT}. 
Let $P$ and $Q$ be linear projections on $Y$ with $\ker P=\ker Q=X$ and 
$\dim \text{\rm Range }P = \dim \text{\rm Range }Q = \dim Y/X$. 
Then 
$$\lim_{n\to\infty}\|T_n P-T_n Q\|_{cb} =0\ .
\tag 3.1$$
Hence 
$$\varlimsup_{n\to \infty} \|T_n P\|_{cb} = \varlimsup_{n\to\infty} 
\|T_nQ\|_{cb}\ .
\tag 3.2$$ 
\endproclaim 

\demo{Proof} 
Let $S_n=T_n Q$ for all $n$. 
Then 
$$T_n -S_n \to 0 \text{ in the SOT.}
\tag 3.3$$ 
Indeed, if $y\in Y$, $(T_n-S_n) (y) = T_n(I-Q)(y)\to0$ in norm since 
$(I-Q)y\in X$. 

Let $F=\text{Range }P$. 
Since $F$ is finite-dimensional, (3.3) yields 
$$\|(T_n-S_n) |F\| \to0 
\tag 3.4$$ 
whence 
$$\|(T_n-S_n) |F\|_{cb} \le (\dim F)\|(T_n -S_n)|F\|\to 0\ .
\tag 3.5$$ 
Now 
$$QP = Q\text{ since } Q(I-P)=0\ . 
\tag 3.6$$ 
Hence 
$$\align
\|T_nP-T_nQ\|_{cb} & = \|T_n P-S_nP\|\ \text{ by (3.6)} \tag 3.7\\
&\le \|(T_n-S_n)|F\|_{cb} \|P\|_{cb} \to 0\ .
\endalign$$
(Note that $P$ is completely bounded since it is a continuous finite 
rank operator.)\qed
\enddemo

\remark{Remark} 
This proof could be given ``more conceptually'' by noting that 
$\ker (P-Q) \supset X$ and hence the operator $T_n(P-Q)$ ``lives'' on 
$Y/X$, a finite-dimensional space; in fact $T_n(P-Q) = T_n((I-Q)-(I-P))\to 0$
in the SOT on $Y/X$, so the $cb$-norms of the sequence $(T_n(P-Q))$, as 
operators on $Y/X$, go to zero. 
\endremark 

The proof of Theorem 3.6 is actually identical to that for Theorem~2.11 
once we draw the following consequence of Lemma~3.7. 

\proclaim{Lemma 3.8} 
Assuming the hypotheses of 3.6, let $X\subset Y_0\subset Y$ with 
$Y_0/X$ finite-dimensional and let $P:Y_0\to Y_0$ be a finite-rank 
projection with $\ker P=X$ and $\text{\rm rank } P= \dim Y_0/X$. 
Then 
$$\varlimsup_{j\to\infty} \|T_jP\|_{cb} \le C\ .
\tag 3.8$$
\endproclaim 

\demo{Proof} 
By hypotheses there exists a projection $Q$ with $\ker P=X$ and 
rank~$Q=\dim Y_0/X$, so that $\|Q\|_{cb} \le C$. 
Hence of course 
$$\varlimsup_{n\to\infty} \|T_nQ\|_{cb} \le C\ .
\tag 3.9$$ 
Now (3.8) follows immediately from Lemma~3.7, in virtue of (3.2).\qed
\enddemo 

\remark{Comment} 
Lemma 3.8 holds for arbitrary $Z_j$'s; i.e., the 
assumption of $\lambda$-injectivity 
is not needed here. 
\endremark 

The proof of Theorem 3.6 is now word for word the same as that for 
Theorem~2.11, except that we replace Corollary~2.13 by Lemma~3.8 in the 
discussion.\qed

We need one last ingredient for the proof of Theorem~3.1 (the main result 
in this section). 

(If $X\subset Y$, we identify $X^{**}$ with $X^{\bot\bot}\subset Y^{**}$.) 

\proclaim{Lemma 3.9}  
Let $X$ and $Y$ be operator spaces with $X\subset Y$, $X^{**}$ isomorphically 
injective, and $Y$ locally reflexive. 
Then $X$ is completely locally complemented in $Y$. 
\endproclaim 

\remark{Remark} 
The proof yields that if $X^{**}$ is $\lambda$-injective and $Y$ is 
$C$-locally reflexive, then for all $\ep>0$, $X$ is completely locally 
$(C+\lambda+1+\ep)$-cocomplemented in $Y$, hence $X$ is completely locally 
$(C+\lambda+2+\ep)$-complemented in $Y$. 
\endremark 

We delay the proof of this lemma, showing instead how we obtain Theorem~3.1. 
In fact, we have the more general

\proclaim{Theorem 3.10} 
Let $Z_1,Z_2,\ldots$ be reflexive $\lambda$-injective operator spaces, 
$Z= (Z_j)_{c_0}$, and $X\subset Y$ be operator spaces with $Y/X$ separable 
and $Y$ locally reflexive. 
Let ${T:X\to Z}$
be a complete surjective isomorphism. 
Then $T$ admits a completely bounded extension $\tilde T:Y\to Z$.
\endproclaim 

\remark{Remarks 1} 
If $Y$ is $C$-locally reflexive and $\|T\|_{cb} \|T^{-1}\|_{cb} =\gamma$, 
we obtain for $\ep>0$ that $\tilde T$ may be chosen with 
$$\|\tilde T\|_{cb} < 
(\gamma\lambda^3 +(C+1)\lambda^2 +\lambda +\ep) 
\|T\|_{cb}\ .$$
\endremark

\remark{2} 
As noted in Section 2, every separable isomorphically injective operator 
space is reflexive, so 3.10 indeed yields Theorem~3.1. 
Actually, more care in the proof yields that 
{\it the conclusion of Theorem~3.10 holds if one deletes the 
reflexivity assumption  from its hypotheses.} 
Hence we obtain the ``quantized'' version of Corollary~1.2 
(with a worse constant): 
{\it Let $Z$ be as in Theorem~3.10 (but drop the assumption that the $Z_j$'s 
are reflexive). 
Then $Z$ is $\lambda^3 +(C+1)\lambda^2 +\lambda +\ep)$-completely 
complemented in every $C$-locally 
reflexive superspace $Y$ with $Y/Z$ separable.} 
\endremark 

\demo{Proof of Theorem~3.10} 
Let $C$ and $\gamma$ be as in Remark~1. 
It follows from the hypotheses that $X^{**}$ is completely 
$\gamma$-isomorphic  to $(Z_j)_{\ell^\infty}$, 
a $\lambda$-injective operator space. 
Hence $X^{**}$ is completely $(\gamma\lambda+1)$-cocomplemented in $Y^{**}$, 
so the proof of Lemma~3.9 yields that given $\ep>0$, $X$ is completely 
locally $(C+\gamma\lambda +1+\ep)$-cocomplemented in $Y$. 
Hence (by playing with $\ep$) we obtain from Theorem~3.4 that the 
extension $\tilde T$ may be chosen with 
$$\displaylines{\hbox{(3.10)}\hfill 
\|\tilde T\|_{cb} \le 
(\gamma\lambda^3 +(C+1)\lambda^2 +\lambda +\ep) 
\|T\|_{cb}\ .
\hfill {\qed}}$$
\enddemo 

We now deal with Lemma 3.9. 
Let us first recall the precise concept of operator space local reflexivity 
(reformulated in the spirit of the original Banach space concept given by 
J.~Lindenstrauss and the author in \cite{LR}, as refined in \cite{JRZ}). 

\definition{Definition 3.2} 
An operator space $X$ is called {\it $C$-locally reflexive\/} if for all 
$\ep>0$, and finite dimensional subspaces $F$ and $G$ of $X^*$ and $X^{**}$ 
respectively, there exists a linear operator $T:G\to X$ satisfying 
$$\langle Tg,f\rangle = \langle g,f\rangle \text{ for all } g\in G\ ,\ 
f\in F
\tag 3.11$$ 
and 
$$\|T\|_{cb} < C+\ep\ .
\tag 3.12$$ 
\enddefinition

As shown in \cite{JRZ}, Banach spaces are thus 1-locally reflexive. 
Remarkable permanence properties given in \cite{ER} yield that if $X$ 
is a $C$-locally reflexive operator space, then every subspace of $X$ is 
$C$-locally reflexive; moreover as noted above, nuclear $C^*$-algebras 
are 1-locally reflexive. 

Lemma 3.8 is an immediate consequence of the following technical result 
(whose proof is the operator space analogue of an argument in \cite{FJT}). 

\proclaim{Sublemma 3.11} 
Let $X\subset Y$ be operator spaces with $\dim Y/X<\infty$  so that $Y$ 
is $C$-locally reflexive and $X^{**}$ is completely $\beta$-cocomplemented 
in $Y^{**}$. 
Then for all $\ep>0$, $X$ is completely $(C+\beta+\ep)$-cocomplemented in $Y$. 
\endproclaim

Let us first deduce Lemma 3.9. 
Let $\ep>0$. 
Assuming that $X^{**}$ is $\lambda$-injective, $X^{**}$ is completely  
$\lambda$-complemented in $Y^{**}$. 
Now assuming $Y$ is $C$-locally reflexive, if $Y_0$ is a subspace of $Y$ 
with $X\subset Y_0$ and $Y_0/X$ finite-dimensional, $Y_0$ is also 
$C$-locally reflexive and of course $X^{**}$ is also completely 
$\lambda$-complemented in $Y_0^{**}$, hence $X^{**}$ is completely 
$(\lambda+1)$-cocomplemented in $Y_0^{**}$. 
Thus by Sublemma~3.11, $X$ is completely $(C+\lambda+1+\ep)$-cocomplemented 
in $Y_0$.\qed 

\demo{Proof of 3.11}
Of course we identify $X^{**}$ with $X^{\bot\bot}$. 
Let $F= X^\bot$. 
The hypotheses actually imply that there exists a projection $P$ from $Y^*$ 
onto $X^\bot$ satisfying 
$$\|P\|_{cb}\le \beta
\tag 3.13$$

Indeed, if $Q$ is a projection on $Y^{**}$ with $\ker Q= X^{**}$ and 
$\|Q\|_{cb}\le\beta$, then $P= Q^*|Y^*$ has the desired property, where we 
regard $Y^* \subset Y^{***}$. 
(In fact, the range of $Q^*$ equals $X^{\bot\bot\bot} = X^\bot$, because 
$Y/X$ is finite-dimensional). 

Now define $G$ by 
$$G = P^* (Y^{**}) \ .
\tag 3.14$$

Of course $G$ is finite-dimensional; hence since $Y$ is $C$-locally reflexive, 
given $\ep>0$, choose $T:G\to Y$ a linear operator with 
$$\|T\|_{cb} < C+ {\ep\over\beta}
\tag 3.15$$
and 
$$\langle Tg,f\rangle = \langle g,f\rangle\ \text{ for all } g\in G\ ,\ 
f\in F\ .
\tag 3.16$$ 
Finally, define $H$ by 
$$H= T(G)\ . 
\tag 3.17$$ 

We now claim that $H$ yields the desired decomposition of $Y$. 
Now it follows immediately from (3.14) that 
$$Y^{**} = F^\bot \oplus G\ .
\tag 3.18$$ 

This and (3.16) imply that $T$ is one-one and $H\cap X= \{0\}$. 
Indeed, suppose $g\in G$ and $Tg=0$. 
Then $\langle Tg,f\rangle = \langle g,f\rangle=0$ for all $f\in F$, whence 
by (3.18), $g=0$. 
But if $Tg\in X$, then since $X^\bot =F$, $\langle Tg,f\rangle =0 = 
\langle g,f\rangle$ for all $f\in F$, so of course $g= Tg=0$. 

Since $\dim Y/X = \dim Y^{**}/X^{**} = \dim G$, we have now deduced 
$$Y = X\oplus H\ . 
\tag 3.19$$
Now let $R$ be the projection from $Y$ onto $H$ with $\ker R=X$. 
We claim 
$$\|R\|_{cb} < C\beta +\ep\ .
\tag 3.20$$ 

We need the fundamental duality pairing for operator spaces. 
Fix $K_1,\ldots,K_m$ in $\bK$. 
Then given $y_1,\ldots,y_m$ in $Y$, $y_1^*,\ldots,y_\ell^*$ in $Y^*$, and 
$L_1,\ldots,L_\ell$ in $\bK$, we define 
$$\left\langle \sum_{i=1}^m K_i\otimes y_i, 
\sum_{j=1}^\ell L_j \otimes y_j^*\right\rangle 
= \sum_{i,j} y_j^* (y_i) K_i \otimes L_j\ . 
\tag 3.21$$ 
(Here, the last term is an operator on $\ell^2 \otimes \ell^2$.) 
Then we have (cf.\ \cite{Pi}) 
$$\|\sum K_i\otimes y_i\| = \sup \left\{ \|\langle \sum K_i\otimes y_i, 
\sum L_j \otimes y_j^*\rangle\| : \|\sum L_j\otimes y_j^*\| =1\right\}\ .
\tag 3.22$$

Now applying this duality statement to $Y^*$ rather than $Y$, it follows 
by our definition of $P$ and $G$, that given $g_1,\ldots,g_m$ in $G$ 
(and $K_1,\ldots, K_m$ as above), then 
$$\align 
\|\sum K_i\otimes g_i\| 
&\le \beta \sup \Big\{ \|\langle \sum K_i\otimes g_i,
\sum L_j\otimes f_j\rangle\| : 
\tag 3.23\\
&\qquad f_1,\ldots, f_\ell \in F,\ L_1,\ldots,L_\ell\text{ in } \bK,\\ 
&\qquad \text{ and } \|\sum L_j\otimes f_j\| = 1\Big\}\ .
\endalign$$
Finally, let $h_1,\ldots,h_n$ in $H$, $x_1,\ldots,x_m$ in $X$ and 
$K_1,\ldots,K_m$ as above. 
We must prove: 
$$\align 
\|\sum K_i\otimes R(x_i+h_i)\| 
& = \|\sum K_i\otimes h_i\|\ \text{ (trivial)}
\tag 3.24\cr 
& \le (C\beta +\ep) \|\sum K_i \otimes (x_i+h_i)\|\ .
\endalign$$

Now choose unique $g_1,\ldots,g_m$ in $G$ with $h_i = Tg_i$ for all $i$. 
Then 
$$\align 
\|\sum K_i \otimes Tg_i\| 
& \le (C+\ep/\beta) \|\sum K_i \otimes g_i\|\qquad\qquad\text{ by (3.15)}\\
\noalign{\vskip-20pt}
&\tag 3.25\\
& \le (C+\ep/\beta) \beta \sup \Big\{\| \langle \sum K_i\otimes g_i,
\sum L_j\otimes f_j\rangle\| : \cr
&\qquad \text{$f_i$'s $\in F$, $L_j$'s in $\bK$, and}\\
&\qquad \|\sum L_j \otimes f_j\| =1\Big\}\qquad\qquad \text{ by (3.23)}\\
& = (C\beta +\ep) \sup \Big\{\| \langle \sum K_i \otimes (x_i+Tg_i), 
\sum L_j \otimes f_j\rangle\| : \cr
&\qquad \text{$f_i$'s $\in F$, $L_j$'s in $\bK$ and } \cr
&\qquad \|\sum L_j\otimes f_j\|=1\Big\}
\endalign$$
by (3.18) and the fact that $F= X^\bot$
$$\le (C\beta+\ep) \|\sum K_i \otimes (x_i+Tg_i)\|\ \text{ by (3.21).}$$ 
This proves (3.24), completing the proof.\qed 
\enddemo

\head 4. Examples of spaces with the CSEP and the CSCP\endhead 

Our preceding results yield lists of separable infinite-dimensional operator 
spaces with the CSEP and CSCP. 
It is conceivable that these lists are complete (up to complete isomorphism). 

The results stated in this section are direct consequences of the work in 
the preceding sections and previously known facts. 
The conjectures we formulate here are strongly believed to be true, and 
should be ``accessible.'' 
On the other hand, the problems we formulate are (probably) at a considerably 
deeper level. 

We first give a basic definition; the operator space analogue of a well 
known Banach space concept. 

\definition{Definition 4.1} 
An operator space $X$ is called {\it primary\/} if whenever $Y$ and $Z$ are 
operator spaces with $X$ completely isomorphic to $Y\oplus Z$, then $X$ is 
completely isomorphic to $Y$ or $Z$. 
\enddefinition 

\noindent 
All of our examples of spaces with the CSEP (resp. CSCP) are direct sums 
of primary spaces with the CSEP (resp. CSCP). 

We first treat the CSEP. 
Recall that $\bR$, $\bC$ denote infinite-dimensional row and column space, 
respectively, and $\bR_n$, $\bC_n$ $n$-dimensional row and column space, 
respectively. 

\proclaim{Proposition 4.1} 
There are at least six isomorphically different Banach spaces among the 
separable infinite-dimensional operator spaces with the {\rm CSEP}, namely  
$$c_0\ ,\quad (\ell_n^2)_{c_0}\ ,\quad  c_0 (\ell^2)\ ,\quad 
\ell^2\ ,\quad c_0\oplus \ell^2\ ,\ \text{ and }\ (\ell_n^2)_{c_0}\oplus 
\ell^2\ .
\tag 4.1$$
\endproclaim 

\demo{Proof} 
Standard Banach space results easily yield these spaces are isomorphically 
distinct (cf.\ \cite{BCLT}). 
Of course Sobczyk's theorem yields that $c_0$ has the 2-CSEP; as a Banach 
space, $R$ is just isometric to $\ell^2$, and of course $R$ has the 1-CSEP. 
Corollary~2.7 yields immediately that $(R_n)_{c_0}$ and $c_0(R)$ have the 
$(2+\ep)$-CSEP for all $\ep>0$, and of course $(R_n)_{c_0}$ is 
isometric to $(\ell_n^2)_{c_0}$ and $c_0(R)$ is isometric to $c_0(\ell^2)$. 
Finally, $c_0\oplus R$ has the 2-CSEP, and of course this is just 
$c_0\oplus \ell^2$ in the Banach space category.\qed
\enddemo 

\proclaim{Problem 4.1} 
Let $X$ be a separable infinite-dimensional operator space with the {\rm CSEP}. 
Is $X$ Banach isomorphic to one of the six spaces in $(4.1)$?
\endproclaim 

By the results in \cite{BCLT}, the first four spaces in (4.1) are primary, 
and moreover every infinite dimensional complemented subspace of $c_0(\ell^2)$ 
(the largest one), is isomorphic to one of these six. 
Thus Problem~4.1 has an affirmative answer if every separable space with 
the CSEP is completely isomorphic to a subspace of $c_0(\bR)\oplus 
c_0(\bC)$. 

\proclaim{Conjecture 4.2} 
There are at least seven completely isomorphically distinct separable 
infinite-dimensional primary operator spaces with the {\rm CSEP}, namely 
$$c_0\ ,\quad (\bR_n)_{c_0}\ ,\quad (\bC_n)_{c_0}\ ,\quad \bR\ ,\quad \bC\ ,
\quad c_0(\bR)\ ,\quad c_0(\bC)\ .
\tag 4.2$$ 
\endproclaim 

As before, it follows immediately from Corollary~2.7 that all these spaces 
have the CSEP (indeed all are completely isometric to completely 
contractively complemented subspaces of $c_0(\bR) \oplus c_0(\bC)$). 
It is also easily seen that all these spaces are isomorphically distinct as 
operator spaces, and it is essentially trivial that $c_0$, $\bR$, and $\bC$ 
are all primary (in fact they are prime). 
The content of the conjecture thus becomes: 
the remaining spaces in (4.2) are all primary. 

\proclaim{Problem 4.2} 
Let $X$ be a separable  infinite-dimensional primary operator space with 
the {\rm CSEP}. 
Is $X$ completely isomorphic to one of the seven spaces listed in $(4.2)$?
\endproclaim 

\proclaim{Conjecture 4.3} 
There are at least $21$ completely 
isomorphically distinct separable operator spaces 
with the {\rm CSEP}, namely 
\roster
\item"(a)" the seven spaces listed in $(4.2)$
\item"(b)" the nine spaces $c_0\oplus \bR$, $c_0\oplus\bC$, $(\bR_n)_{c_0} 
\oplus \bR$, $(\bR_n)_{c_0} \oplus (\bC_n)_{c_0}$, 
$(\bR_n)_{c_0}\oplus \bC$, 
$(\bC_n)_{c_0} \oplus \bR$, 
$\bR\oplus \bC$, 
$(\bC_n)_{c_0} \oplus \bC$, 
$c_0(\bR)\oplus c_0(\bC)$ 
\item"(c)" the five spaces $c_0\oplus \bR\oplus \bC$, 
$(\bR_n)_{c_0} \oplus (\bC_n)_{c_0} \oplus \bR$, 
$(\bR_n)_{c_0} \oplus (\bC_n)_{c_0} \oplus \bC$, 
$(\bR_n)_{c_0} \oplus \bR\oplus \bC$, 
$(\bR_n)_{c_0} \oplus (\bC_n)_{c_0} \oplus \bR\oplus \bC$.
\endroster 
Moreover any finite direct sum of any of these spaces is again completely 
isomorphic to one of them. 
\endproclaim 

As above, it follows immediately from the results of Section~2 that all these 
spaces have the CSEP. 
We leave the remaining assertions of this conjecture to the ambitious reader. 

\proclaim{Problem 4.3} 
Is every separable infinite-dimensional operator space with the {\rm CSEP} 
completely isomorphic to one of the $21$ spaces in Conjecture~4.3?
\endproclaim 

We now deal with the CSCP. 
It is conceivable that the separable infinite-dimensional operator spaces 
with the CSCP are precisely those which are completely isomorphic  to a 
completely complemented subspaces of $\bK$. 
Accordingly, we discuss the evident spaces with this property; recall that 
$\bK_0$ denotes the space $(M_n)_{c_0}$. 
The following result is due to J.~Arazy and  J.~Lindenstrauss  (see 
Theorem~5 and Remark (i), p.107, of \cite{AL}). 

\proclaim{Proposition 4.4} 
There are at least $11$ isomorphically distinct Banach spaces isomorphic 
to an infinite-dimensional complemented subspace of $\bK$, namely 
\roster
\item"(a)" the seven spaces  
$c_0$, $\ell^2$, $(\ell_n^2)_{c_0}$, $c_0(\ell^2)$, $\bK_0$,  
$(M_{\infty,n})_{c_0}$, and $\bK$ 
\item"(b)" the four spaces $c_0\oplus \ell^2$, 
$\ell^2 \oplus (\ell_n^2)_{c_0}$,$\ell^2 \oplus \bK_0$,  and 
$c_0 (\ell^2)\oplus \bK_0$. 
\endroster
\endproclaim 

It is known that all the spaces in (a), except possibly $(M_{\infty,n})_{c_0}$, 
are primary. 
The primariness of the first four is noted above (\cite{BCLT}). 
The result that $\bK$ and $\bK_0$ are primary, is due to 
J.~Arazy \cite{Ar}. 
I conjecture that also $(M_{\infty,n})_{c_0}$ is primary, but this remains 
an open question. 

\proclaim{Problem 4.4} 
Is every infinite-dimensional completely complemented subspace of $\bK$ 
Banach-isomorphic to one of the $11$ spaces listed in 
{\rm (a)} and {\rm (b)} of $4.4$?
\endproclaim 

It is conceivable that every infinite dimensional complemented subspace of 
$\bK$ is isomorphic to one of these 11 spaces; this problem is raised 
in \cite{AL}. 
Problem~4.4 might be somewhat more accessible. 
Of course our motivation here is that by the results of Section~3 
(resp. \cite{OR} for $\bK$ itself), all of the spaces listed in 4.4  are 
Banach-isomorphic to operator spaces with the CSCP. 

\proclaim{Conjecture 4.5} 
There are at least $11$ completely 
isomorphically distinct primary operator spaces, each 
completely isometric to a completely contractively complemented subspace  
of $\bK$, namely 
\roster
\item"(a)" the seven spaces listed in {\rm (4.2)} 
\item"(b)" the four spaces $\bK$, $\bK_0$, $(M_{\infty,n})_{c_0}$ and 
$(M_{n,\infty})_{c_0}$. 
\endroster
\endproclaim 

Using the known Banach space result, Proposition~4.4, 
it is not hard to see that all the listed spaces are 
completely isomorphically distinct, 
and all are completely contractively complemented in $\bK$. 
The content of the conjecture thus becomes: 
all these spaces are primary. 
(It seems likely the work in \cite{Ar} should yield that $\bK$ and $\bK_0$ 
are primary operator spaces, but we have not verified this.) 
Again, by the results of Section~3, 
(and \cite{OR} for the case of $\bK$ itself) 
all these spaces  have the CSCP. 

\proclaim{Problem 4.5} 
Let $X$ be a separable infinite-dimensional primary operator space with 
the {\rm CSCP}. 
Is $X$ completely isomorphic to one of the spaces listed in {\rm (a)} 
and {\rm (b)} of Conjecture~4.5?
\endproclaim 

Of course a motivation to classify the (apparently finite but rather 
immense number of) finite-direct sums of these 11 spaces would be provided 
by an affirmative answer to the following (obviously deep) problem: 

\proclaim{Problem 4.6} 
Is every operator space with the {\rm CSCP} completely isomorphic to a 
finite direct sum of primary operator spaces?
\endproclaim 

\subhead Acknowledgments\endsubhead 
The research for this paper was partially supported by NSF Grant 
DMS-9500874  and TARP Grant ARP-275.
\enddocument

\vskip.2in 
{\parskip=0pt 
\baselineskip=12pt
\item{} Haskell Rosenthal
\item{} {\sl Department of Mathematics}
\item{} {\sl The University of Texas at Austin}
\item{} {\sl Austin, TX 78712}
\item{} {\sl E-mail\/}: {\tt  rosenthl\@math.utexas.edu}
\smallskip
}


\Refs 
\widestnumber\key{BCLT}

\ref\key A
\by T. Andersen 
\paper Linear extensions, projections, and split faces 
\jour J. Funct. Anal. \vol 17 \yr 1974 \pages 161--173 
\endref 

\ref\key An  
\by T. Ando 
\paper Closed range theorems for convex sets and linear liftings 
\jour Pacific J. Math. \vol 44 \yr 1973 \pages 393--410 
\endref 

\ref\key Ar 
\by J. Arazy 
\paper A remark on complemented subspaces of unitary matrix spaces 
\jour Proc. Amer. Math. Soc. \vol 79 \yr 1980 \pages 601--608 
\endref 

\ref\key AL 
\by J. Arazy and J. Lindenstrauss 
\paper Some linear topological properties of the spaces $C_p$ of operators 
on Hilbert space 
\jour Compositio Math. \vol 30 \yr 1975 \pages 81--111 
\endref 

\ref\key Arv 
\by W.B. Arveson 
\paper Subalgebras of $C^*$-algebras 
\jour Acta Math. \vol 123 \yr 1969 \pages 141--224 
\endref

\ref\key BP 
\by D. Blecher and V. Paulsen 
\paper Tensor products of operator spaces 
\jour J. Funct. Anal. \vol 99 \yr 1991 \pages 262--292 
\endref 

\ref\key BCLT 
\by J. Bourgain, P.G. Casazza, J. Lindenstrauss and L. Tzafriri 
\paper Banach spaces with a unique unconditional basis up to permutation 
\jour Memoirs Amer. Math. Soc. 
\finalinfo No. 322, 1985 
\endref 

\ref\key CE 
\by M.-D. Choi and E. Effros 
\paper The completely positive lifting problem for $C^*$-algebras 
\jour Ann. of Math. \vol 104 \yr 1976 \pages 585--609
\endref 

\ref\key EH 
\by E. Effros and U. Haagerup 
\paper Lifting problems and local reflexivity for $C^*$-algebras 
\jour Duke Math. J. \vol 52 \yr 1985 \pages 103--128 
\endref 

\ref\key ER1 
\by E. Effros and Z.J. Ruan 
\paper On the abstract characterization of operator spaces 
\jour Proc. Amer. Math. Soc. \vol 119 \yr 1993 \pages 579--584 
\endref 

\ref\key ER2 
\by E. Effros and Z.J. Ruan 
\paper Mapping spaces and liftings for operator spaces 
\jour Proc. London Math. Soc. \vol 69 \yr 1994 \pages 171--197 
\endref 

\ref\key FJT 
\by T. Figiel, W. Johnson and L.Tzafriri
\paper On Banach lattices and spaces having local unconditional 
structure, with applications to Lorentz function spaces 
\jour J. Approx. Theory \vol 13 \yr 1975 \pages 27--48 
\endref 

\ref\key HWW
\by P. Harmand, D. Werner and W. Werner  
\book $M$-ideals in Banach spaces and Banach algebras
\bookinfo SLNM 1547 
\publ Springer-Verlag 
\yr 1993 
\endref 

\ref\key J 
\by R.C. James 
\paper Uniformly non-square Banach spaces 
\jour Annals of Math. \vol 80 \yr 1964 \pages 542--550 
\endref 

\ref\key JO 
\by W. Johnson and T. Oikhberg 
\paper Relative complementation property 
\finalinfo preprint
\endref 

\ref\key JRZ 
\by W.B. Johnson, H.P. Rosenthal and M. Zippin 
\paper On bases, finite dimensional decompositions and weaker structures 
in Banach spaces 
\jour Israel J. Math. \vol 9 \yr 1971 \pages 488--506 
\endref 

\ref\key K 
\by E. Kirchberg 
\paper On non-semisplit extensions, tensor products and exactness of 
group $C^*$-algebras 
\jour Invent. Math. \vol 112 \yr 1993 \pages 449--489
\endref 

\ref\key LR 
\by J. Lindenstrauss and H.P. Rosenthal 
\paper The ${\Cal L}_p$ spaces 
\jour Israel J. Math. \vol 7 \yr 1969 \pages 325--349 
\endref 

\ref\key M 
\by A. Martinez-Abejon 
\paper A short proof of the principle of local reflexivity 
\jour Proc. Amer. Math. Soc. 
\toappear 
\endref 

%\ref\key O1 
%\by T. Oikhberg 
%\paper Projections onto some spaces of sequences 
%\finalinfo preprint 
%\endref 

\ref\key O 
\by T. Oikhberg 
\paper Direct sums of operator spaces 
\finalinfo preprint 
\endref 

\ref\key OR
\by T. Oikhberg and H.P. Rosenthal 
\paper Some extension properties for the space of compact operators 
\finalinfo in preparation 
\endref

\ref\key P 
\by  V. Paulsen 
\book Completely Bounded Maps and Dilations 
\publ John Wiley \& Sons, Inc. 
\publaddr New York 
\yr 1986
\endref 

\ref\key Pe 
\by A. Pe{\l}czy\'nski 
\paper Projections in certain Banach spaces 
\jour Studia Math. \vol 29 \yr 1969 \pages 209--227 
\endref 

\ref\key Pf
\by H. Pfitzner 
\paper Weak compactness in the dual of a $C^*$-algebra is determined 
commutatively 
\jour Math. Ann. \vol 298 \pages 349--371 \yr 1994 
\endref 

\ref\key Pi 
\by G. Pisier 
\paper An introduction to the theory of operator spaces 
\finalinfo preprint 
\endref 

\ref\key R 
\by A. Guyan Robertson 
\paper Injective matricial Hilbert spaces 
\jour Math. Proc. Camb. Phil. Soc. \vol 110 \yr 1991 \pages  183--190 
\endref 

\ref\key Ro 
\by H.P. Rosenthal 
\paper On complemented and quasi-complemented subspaces of quotients of 
$C(S)$ for Stonian $S$ 
\jour Proc. Natl. Acad. Sci. \vol 60 \yr 1968 \pages 1165--1169
\endref 

\ref\key Ru 
\by Z.J. Ruan 
\paper Injectivity and operator spaces 
\jour Trans. Amer. Math. Soc. \vol 315 \yr 1989 \pages 89--104 
\endref 

\ref\key S 
\by A. Sobczyk 
\paper Projection of the space $(m)$ on its subspace $(c_0)$ 
\jour Bull. Amer. Math. Soc. \vol 47 \yr 1941 \pages 938--947 
\finalinfo MR3-205
\endref 

\ref\key S 
\by R.R. Smith 
\paper Completely bounded mpas between $C^*$-algebras 
\jour J. London Math. Soc (2) \vol 27 \yr 1983 \pages 157--166 
\endref 

\ref\key SW 
\by R.R. Smith and D.P. Williams 
\paper Separable Injectivity for $C^*$-algebras 
\jour Indiana Univ. Math. J. 
\vol 37 \yr 1988 \pages 111--133 
\endref 

\ref\key V 
\by W.A. Veech 
\paper Short proof of Sobczyk's theorem 
\jour Proc. Amer. Math. Soc. \vol 28 \yr 1971 \page 627--628
\endref 

\ref\key W 
\by S. Wassermann 
\paper Exact $C^*$-algebras and related topics 
\jour Seoul National University Lecture Notes Series 
\finalinfo No. 19, 1994 
\endref 

\ref\key Wi 
\by G. Wittstock 
\paper Ein operatorwertigen Hahn-Banach Satz 
\jour J. Func. Anal. \vol 40 \yr 1981 \pages 127--150 
\endref 

\ref\key Z 
\by M. Zippin 
\paper The separable extension problem 
\jour Israel J. Math. \vol 26 \yr 1977 \pages 372--387 
\endref 

\endRefs
\enddocument